
\catcode'32=9
\magnification=1200

\voffset=1cm

\font\tenpc=cmcsc10

\font\eightrm=cmr8
\font\eighti=cmmi8
\font\eightsy=cmsy8
\font\eightbf=cmbx8
\font\eighttt=cmtt8
\font\eightit=cmti8
\font\eightsl=cmsl8
\font\sixrm=cmr6
\font\sixi=cmmi6
\font\sixsy=cmsy6
\font\sixbf=cmbx6

\skewchar\eighti='177 \skewchar\sixi='177
\skewchar\eightsy='60 \skewchar\sixsy='60

\font\tengoth=eufm10
\font\tenbboard=msbm10
\font\eightgoth=eufm7 at 8pt
\font\eightbboard=msbm7 at 8pt
\font\sevengoth=eufm7
\font\sevenbboard=msbm7
\font\sixgoth=eufm5 at 6 pt
\font\fivegoth=eufm5

\font\tengoth=eufm10
\font\tenbboard=msbm10
\font\eightgoth=eufm7 at 8pt
\font\eightbboard=msbm7 at 8pt
\font\sevengoth=eufm7
\font\sevenbboard=msbm7
\font\sixgoth=eufm5 at 6 pt
\font\fivegoth=eufm5

\newfam\gothfam
\newfam\bboardfam

\catcode`\@=11

\def\raggedbottom{\topskip 10pt plus 36pt
\r@ggedbottomtrue}
\def\pc#1#2|{{\bigf@ntpc #1\penalty
\@MM\hskip\z@skip\smallf@ntpc #2}}

\def\tenpoint{%
  \textfont0=\tenrm \scriptfont0=\sevenrm \scriptscriptfont0=\fiverm
  \def\rm{\fam\z@\tenrm}%
  \textfont1=\teni \scriptfont1=\seveni \scriptscriptfont1=\fivei
  \def\oldstyle{\fam\@ne\teni}%
  \textfont2=\tensy \scriptfont2=\sevensy \scriptscriptfont2=\fivesy
  \textfont\gothfam=\tengoth \scriptfont\gothfam=\sevengoth
  \scriptscriptfont\gothfam=\fivegoth
  \def\goth{\fam\gothfam\tengoth}%
  \textfont\bboardfam=\tenbboard \scriptfont\bboardfam=\sevenbboard
  \scriptscriptfont\bboardfam=\sevenbboard
  \def\bboard{\fam\bboardfam}%
  \textfont\itfam=\tenit
  \def\it{\fam\itfam\tenit}%
  \textfont\slfam=\tensl
  \def\sl{\fam\slfam\tensl}%
  \textfont\bffam=\tenbf \scriptfont\bffam=\sevenbf
  \scriptscriptfont\bffam=\fivebf
  \def\bf{\fam\bffam\tenbf}%
  \textfont\ttfam=\tentt
  \def\tt{\fam\ttfam\tentt}%
  \abovedisplayskip=12pt plus 3pt minus 9pt
  \abovedisplayshortskip=0pt plus 3pt
  \belowdisplayskip=12pt plus 3pt minus 9pt
  \belowdisplayshortskip=7pt plus 3pt minus 4pt
  \smallskipamount=3pt plus 1pt minus 1pt
  \medskipamount=6pt plus 2pt minus 2pt
  \bigskipamount=12pt plus 4pt minus 4pt
  \normalbaselineskip=12pt
  \setbox\strutbox=\hbox{\vrule height8.5pt depth3.5pt width0pt}%
  \let\bigf@ntpc=\tenrm \let\smallf@ntpc=\sevenrm
  \let\petcap=\tenpc
  \normalbaselines\rm}
\def\eightpoint{%
  \textfont0=\eightrm \scriptfont0=\sixrm \scriptscriptfont0=\fiverm
  \def\rm{\fam\z@\eightrm}%
  \textfont1=\eighti \scriptfont1=\sixi \scriptscriptfont1=\fivei
  \def\oldstyle{\fam\@ne\eighti}%
  \textfont2=\eightsy \scriptfont2=\sixsy \scriptscriptfont2=\fivesy
  \textfont\gothfam=\eightgoth \scriptfont\gothfam=\sixgoth
  \scriptscriptfont\gothfam=\fivegoth
  \def\goth{\fam\gothfam\eightgoth}%
  \textfont\bboardfam=\eightbboard \scriptfont\bboardfam=\sevenbboard
  \scriptscriptfont\bboardfam=\sevenbboard
  \def\bboard{\fam\bboardfam}%
  \textfont\itfam=\eightit
  \def\it{\fam\itfam\eightit}%
  \textfont\slfam=\eightsl
  \def\sl{\fam\slfam\eightsl}%
  \textfont\bffam=\eightbf \scriptfont\bffam=\sixbf
  \scriptscriptfont\bffam=\fivebf
  \def\bf{\fam\bffam\eightbf}%
  \textfont\ttfam=\eighttt
  \def\tt{\fam\ttfam\eighttt}%
  \abovedisplayskip=9pt plus 2pt minus 6pt
  \abovedisplayshortskip=0pt plus 2pt
  \belowdisplayskip=9pt plus 2pt minus 6pt
  \belowdisplayshortskip=5pt plus 2pt minus 3pt
  \smallskipamount=2pt plus 1pt minus 1pt
  \medskipamount=4pt plus 2pt minus 1pt
  \bigskipamount=9pt plus 3pt minus 3pt
  \normalbaselineskip=9pt
  \setbox\strutbox=\hbox{\vrule height7pt depth2pt width0pt}%
  \let\bigf@ntpc=\eightrm \let\smallf@ntpc=\sixrm
  \normalbaselines\rm}

\tenpoint

\frenchspacing


\newif\ifpagetitre
\newtoks\auteurcourant \auteurcourant={\hfil}
\newtoks\titrecourant \titrecourant={\hfil}

\def\appeln@te{}
\def\vfootnote#1{\def\@parameter{#1}\insert\footins\bgroup\eightpoint
  \interlinepenalty\interfootnotelinepenalty
  \splittopskip\ht\strutbox 
  \splitmaxdepth\dp\strutbox \floatingpenalty\@MM
  \leftskip\z@skip \rightskip\z@skip
  \ifx\appeln@te\@parameter\indent \else{\noindent #1\ }\fi
  \footstrut\futurelet\next\fo@t}

\pretolerance=500 \tolerance=1000 \brokenpenalty=5000
\newdimen\hmargehaute \hmargehaute=0cm
\newdimen\lpage \lpage=13.3cm
\newdimen\hpage \hpage=20cm
\newdimen\lmargeext \lmargeext=1cm
\hsize=11.25cm
\vsize=18cm
\parskip 0pt
\parindent=12pt

\def\margehaute{\vbox to \hmargehaute{\vss}}%
\def\margebasse{\vss}

\output{\shipout\vbox to \hpage{\margehaute\nointerlineskip
  \corpsdepage\margebasse}
  \advancepageno \global\pagetitrefalse
  \ifnum\outputpenalty>-20000 \else\dosupereject\fi}

\def\corpsdepage{\hbox to \lpage{\hss\pagetexte\hskip\lmargeext}}
\def\pagetexte{\vbox{\makeheadline\pagebody\makefootline}}
\headline={\ifpagetitre\titleheadline \else
  \ifodd\pageno\rightheadline \else\leftheadline\fi\fi}
\def\leftheadline{\eightpoint\hfil\the\auteurcourant\hfil}
\def\rightheadline{\eightpoint\hfil\the\titrecourant\hfil}
\def\titleheadline{\hfill}
\pagetitretrue

\def\footnoterule{\kern-6\p@
  \hrule width 2truein \kern 5.6\p@} 

\def\pd#1#2 {\pc#1#2| }

\def\pointir{\discretionary{.}{}{.\kern.35em---\kern.7em}\nobreak
\hskip 0em plus .3em minus .4em }

\def\abstract#1{\vbox{\eightpoint \pc ABSTRACT|\pointir #1}}

\def\titre#1|{\message{#1}
              \par\vskip 30pt plus 24pt minus 3pt\penalty -1000
              \vskip 0pt plus -24pt minus 3pt\penalty -1000
              \centerline{\bf #1}
              \vskip 5pt
              \penalty 10000 }

\def\section#1|{\par\vskip .3cm
                {\bf #1}\pointir}

\def\ssection#1|{\par\vskip .2cm
                {\it #1}\pointir}

\long\def\th#1|#2\finth{\par\medskip
              {\petcap #1\pointir}{\it #2}\par\smallskip}

\long\def\tha#1|#2\fintha{\par\medskip
                    {\petcap #1.}\par\nobreak{\it #2}\par\smallskip}
\def\cf{{\it cf}}

\def\rem#1|{\par\medskip
            {{\it #1}.\quad}}

\def\rema#1|{\par\medskip
             {{\it #1.}\par\nobreak }}

\def\article#1|#2|#3|#4|#5|#6|#7|
    {{\leftskip=7mm\noindent
     \hangindent=2mm\hangafter=1
     \llap{[#1]\hskip.35em}{#2}.\quad
     #3, {\sl #4}, vol.\nobreak\ {\bf #5}, {\oldstyle #6},
     p.\nobreak\ #7.\par}}
\def\livre#1|#2|#3|#4|
    {{\leftskip=7mm\noindent
    \hangindent=2mm\hangafter=1
    \llap{[#1]\hskip.35em}{#2}.\quad
    {\sl #3}.\quad #4.\par}}
\def\divers#1|#2|#3|
    {{\leftskip=7mm\noindent
    \hangindent=2mm\hangafter=1
     \llap{[#1]\hskip.35em}{#2}.\quad
     #3.\par}}
\mathchardef\conj="0365
\def\proof{\par{\it Proof}.\quad}
\def\qed{\quad\raise -2pt\hbox{\vrule\vbox to 10pt{\hrule width 4pt
\vfill\hrule}\vrule}}

\def\cqfd{\penalty 500 \hbox{\qed}\par\smallskip}
\def\decale#1|{\par\noindent\hskip 28pt\llap{#1}\kern 5pt}

\catcode`\@=12


\catcode`\@=11
\def\matrice#1{\null \,\vcenter {\normalbaselines \m@th
\ialign {\hfil $##$\hfil &&\  \hfil $##$\hfil\crcr
\mathstrut \crcr \noalign {\kern -\baselineskip } #1\crcr
\mathstrut \crcr \noalign {\kern -\baselineskip }}}\,}

\def\petitematrice#1{\left(\null\vcenter {\normalbaselines \m@th
\ialign {\hfil $##$\hfil 
&&\thinspace  \hfil $##$\hfil\crcr
\mathstrut \crcr \noalign {\kern -\baselineskip } #1\crcr
\mathstrut \crcr \noalign {\kern -\baselineskip }}}\right)}

\catcode`\@=12

\def\qed{\quad\raise -2pt\hbox{\vrule\vbox to 10pt{\hrule width 4pt
   \vfill\hrule}\vrule}}

\def\cqfd{\penalty 500 \hbox{\qed}\par\smallskip}


%


\def\il{\bigl]\kern-.25em\bigl]}
\def\ir{\bigr]\kern-.25em\bigr]}

\def\iil{\bigl>\kern-.25em\bigl>}
\def\iir{\bigr>\kern-.25em\bigr>}



\def\Grille{\setbox13=\vbox to 5mm{\hrule width 110mm\vfill}
\setbox13=\vbox{\offinterlineskip
   \copy13\copy13\copy13\copy13\copy13\copy13\copy13\copy13
   \copy13\copy13\copy13\copy13\box13\hrule width 110mm}
\setbox14=\hbox to 5mm{\vrule height 65mm\hfill}
\setbox14=\hbox{\copy14\copy14\copy14\copy14\copy14\copy14
   \copy14\copy14\copy14\copy14\copy14\copy14\copy14\copy14
   \copy14\copy14\copy14\copy14\copy14\copy14\copy14\copy14\box14}
\ht14=0pt\dp14=0pt\wd14=0pt
\setbox13=\vbox to 0pt{\vss\box13\offinterlineskip\box14}
\wd13=0pt\box13}


\def\fleche(#1,#2)\dir(#3,#4)\long#5{%
\noalign{\nointerlineskip\leftput(#1,#2){\vector(#3,#4){#5}}\nointerlineskip}}


\def\hfl#1#2#3{\smash{\mathop{\hbox to#3{\rightarrowfill}}\limits
^{\scriptstyle#1}_{\scriptstyle#2}}}

\def\gfl#1#2#3{\smash{\mathop{\hbox to#3{\leftarrowfill}}\limits
^{\scriptstyle#1}_{\scriptstyle#2}}}


 \message{`lline' & `vector' macros from LaTeX}
 \catcode`@=11
\def\{{\relax\ifmmode\lbrace\else$\lbrace$\fi}
\def\}{\relax\ifmmode\rbrace\else$\rbrace$\fi}
\def\newcount{\alloc@0\count\countdef\insc@unt}
\def\newdimen{\alloc@1\dimen\dimendef\insc@unt}
\def\newwrite{\alloc@7\write\chardef\sixt@@n}

\newwrite\@unused
\newcount\@tempcnta
\newcount\@tempcntb
\newdimen\@tempdima
\newdimen\@tempdimb
\newbox\@tempboxa

\def\@spaces{\space\space\space\space}
\def\@whilenoop#1{}
\def\@whiledim#1\do #2{\ifdim #1\relax#2\@iwhiledim{#1\relax#2}\fi}
\def\@iwhiledim#1{\ifdim #1\let\@nextwhile=\@iwhiledim
        \else\let\@nextwhile=\@whilenoop\fi\@nextwhile{#1}}
\def\@badlinearg{\@latexerr{Bad \string\line\space or \string\vector
   \space argument}}
\def\@latexerr#1#2{\begingroup
\edef\@tempc{#2}\expandafter\errhelp\expandafter{\@tempc}%
\def\@eha{Your command was ignored.
^^JType \space I <command> <return> \space to replace it
  with another command,^^Jor \space <return> \space to continue without it.}
\def\@ehb{You've lost some text. \space \@ehc}
\def\@ehc{Try typing \space <return>
  \space to proceed.^^JIf that doesn't work, type \space X <return> \space to
  quit.}
\def\@ehd{You're in trouble here.  \space\@ehc}

\typeout{LaTeX error. \space See LaTeX manual for explanation.^^J
 \space\@spaces\@spaces\@spaces Type \space H <return> \space for
 immediate help.}\errmessage{#1}\endgroup}
\def\typeout#1{{\let\protect\string\immediate\write\@unused{#1}}}

\font\tenln    = line10
\font\tenlnw   = linew10

\newdimen\@wholewidth
\newdimen\@halfwidth
\newdimen\unitlength 

\unitlength =1pt


\def\thinlines{\let\@linefnt\tenln \let\@circlefnt\tencirc
  \@wholewidth\fontdimen8\tenln \@halfwidth .5\@wholewidth}
\def\thicklines{\let\@linefnt\tenlnw \let\@circlefnt\tencircw
  \@wholewidth\fontdimen8\tenlnw \@halfwidth .5\@wholewidth}

\def\linethickness#1{\@wholewidth #1\relax \@halfwidth .5\@wholewidth}

\newif\if@negarg

\def\lline(#1,#2)#3{\@xarg #1\relax \@yarg #2\relax
\@linelen=#3\unitlength
\ifnum\@xarg =0 \@vline
  \else \ifnum\@yarg =0 \@hline \else \@sline\fi
\fi}

\def\@sline{\ifnum\@xarg< 0 \@negargtrue \@xarg -\@xarg \@yyarg -\@yarg
  \else \@negargfalse \@yyarg \@yarg \fi
\ifnum \@yyarg >0 \@tempcnta\@yyarg \else \@tempcnta -\@yyarg \fi
\ifnum\@tempcnta>6 \@badlinearg\@tempcnta0 \fi
\setbox\@linechar\hbox{\@linefnt\@getlinechar(\@xarg,\@yyarg)}%
\ifnum \@yarg >0 \let\@upordown\raise \@clnht\z@
   \else\let\@upordown\lower \@clnht \ht\@linechar\fi
\@clnwd=\wd\@linechar
\if@negarg \hskip -\wd\@linechar \def\@tempa{\hskip -2\wd\@linechar}\else
     \let\@tempa\relax \fi
\@whiledim \@clnwd <\@linelen \do
  {\@upordown\@clnht\copy\@linechar
   \@tempa
   \advance\@clnht \ht\@linechar
   \advance\@clnwd \wd\@linechar}%
\advance\@clnht -\ht\@linechar
\advance\@clnwd -\wd\@linechar
\@tempdima\@linelen\advance\@tempdima -\@clnwd
\@tempdimb\@tempdima\advance\@tempdimb -\wd\@linechar
\if@negarg \hskip -\@tempdimb \else \hskip \@tempdimb \fi
\multiply\@tempdima \@m
\@tempcnta \@tempdima \@tempdima \wd\@linechar \divide\@tempcnta \@tempdima
\@tempdima \ht\@linechar \multiply\@tempdima \@tempcnta
\divide\@tempdima \@m
\advance\@clnht \@tempdima
\ifdim \@linelen <\wd\@linechar
   \hskip \wd\@linechar
  \else\@upordown\@clnht\copy\@linechar\fi}

\def\@hline{\ifnum \@xarg <0 \hskip -\@linelen \fi
\vrule height \@halfwidth depth \@halfwidth width \@linelen
\ifnum \@xarg <0 \hskip -\@linelen \fi}

\def\@getlinechar(#1,#2){\@tempcnta#1\relax\multiply\@tempcnta 8
\advance\@tempcnta -9 \ifnum #2>0 \advance\@tempcnta #2\relax\else
\advance\@tempcnta -#2\relax\advance\@tempcnta 64 \fi
\char\@tempcnta}

\def\vector(#1,#2)#3{\@xarg #1\relax \@yarg #2\relax
\@linelen=#3\unitlength
\ifnum\@xarg =0 \@vvector
  \else \ifnum\@yarg =0 \@hvector \else \@svector\fi
\fi}

\def\@hvector{\@hline\hbox to 0pt{\@linefnt
\ifnum \@xarg <0 \@getlarrow(1,0)\hss\else
    \hss\@getrarrow(1,0)\fi}}

\def\@vvector{\ifnum \@yarg <0 \@downvector \else \@upvector \fi}

\def\@svector{\@sline
\@tempcnta\@yarg \ifnum\@tempcnta <0 \@tempcnta=-\@tempcnta\fi
\ifnum\@tempcnta <5
  \hskip -\wd\@linechar
  \@upordown\@clnht \hbox{\@linefnt  \if@negarg
  \@getlarrow(\@xarg,\@yyarg) \else \@getrarrow(\@xarg,\@yyarg) \fi}%
\else\@badlinearg\fi}

\def\@getlarrow(#1,#2){\ifnum #2 =\z@ \@tempcnta='33\else
\@tempcnta=#1\relax\multiply\@tempcnta \sixt@@n \advance\@tempcnta
-9 \@tempcntb=#2\relax\multiply\@tempcntb \tw@
\ifnum \@tempcntb >0 \advance\@tempcnta \@tempcntb\relax
\else\advance\@tempcnta -\@tempcntb\advance\@tempcnta 64
\fi\fi\char\@tempcnta}

\def\@getrarrow(#1,#2){\@tempcntb=#2\relax
\ifnum\@tempcntb < 0 \@tempcntb=-\@tempcntb\relax\fi
\ifcase \@tempcntb\relax \@tempcnta='55 \or
\ifnum #1<3 \@tempcnta=#1\relax\multiply\@tempcnta
24 \advance\@tempcnta -6 \else \ifnum #1=3 \@tempcnta=49
\else\@tempcnta=58 \fi\fi\or
\ifnum #1<3 \@tempcnta=#1\relax\multiply\@tempcnta
24 \advance\@tempcnta -3 \else \@tempcnta=51\fi\or
\@tempcnta=#1\relax\multiply\@tempcnta
\sixt@@n \advance\@tempcnta -\tw@ \else
\@tempcnta=#1\relax\multiply\@tempcnta
\sixt@@n \advance\@tempcnta 7 \fi\ifnum #2<0 \advance\@tempcnta 64 \fi
\char\@tempcnta}

\def\@vline{\ifnum \@yarg <0 \@downline \else \@upline\fi}

\def\@upline{\hbox to \z@{\hskip -\@halfwidth \vrule
  width \@wholewidth height \@linelen depth \z@\hss}}

\def\@downline{\hbox to \z@{\hskip -\@halfwidth \vrule
  width \@wholewidth height \z@ depth \@linelen \hss}}

\def\@upvector{\@upline\setbox\@tempboxa\hbox{\@linefnt\char'66}\raise
     \@linelen \hbox to\z@{\lower \ht\@tempboxa\box\@tempboxa\hss}}

\def\@downvector{\@downline\lower \@linelen
      \hbox to \z@{\@linefnt\char'77\hss}}

\thinlines

\newcount\@xarg
\newcount\@yarg
\newcount\@yyarg
\newcount\@multicnt
\newdimen\@xdim
\newdimen\@ydim
\newbox\@linechar
\newdimen\@linelen
\newdimen\@clnwd
\newdimen\@clnht
\newdimen\@dashdim
\newbox\@dashbox
\newcount\@dashcnt
 \catcode`@=12


\newbox\tbox
\newbox\tboxa

\def\leftzer#1{\setbox\tbox=\hbox to 0pt{#1\hss}%
     \ht\tbox=0pt \dp\tbox=0pt \box\tbox}

\def\rightzer#1{\setbox\tbox=\hbox to 0pt{\hss #1}%
     \ht\tbox=0pt \dp\tbox=0pt \box\tbox}

\def\centerzer#1{\setbox\tbox=\hbox to 0pt{\hss #1\hss}%
     \ht\tbox=0pt \dp\tbox=0pt \box\tbox}

%
\def\image(#1,#2)#3{\vbox to #1{\offinterlineskip
    \vss #3 \vskip #2}}


\def\leftput(#1,#2)#3{\setbox\tboxa=\hbox{%
    \kern #1\unitlength
    \raise #2\unitlength\hbox{\leftzer{#3}}}%
    \ht\tboxa=0pt \wd\tboxa=0pt \dp\tboxa=0pt\box\tboxa}

\def\rightput(#1,#2)#3{\setbox\tboxa=\hbox{%
    \kern #1\unitlength
    \raise #2\unitlength\hbox{\rightzer{#3}}}%
    \ht\tboxa=0pt \wd\tboxa=0pt \dp\tboxa=0pt\box\tboxa}

\def\centerput(#1,#2)#3{\setbox\tboxa=\hbox{%
    \kern #1\unitlength
    \raise #2\unitlength\hbox{\centerzer{#3}}}%
    \ht\tboxa=0pt \wd\tboxa=0pt \dp\tboxa=0pt\box\tboxa}

\unitlength=1mm

\def\put(#1,#2)#3{\noalign{\nointerlineskip
                               \centerput(#1,#2){$#3$}
                                \nointerlineskip}}
\def\segment(#1,#2)\dir(#3,#4)\long#5{%
\leftput(#1,#2){\lline(#3,#4){#5}}}
\def\Deltaa{\mathop{\hbox{$\Delta$}}\limits}
\def\brullet{{\scriptscriptstyle\bullet}}

\def\Eoc{\mathop{\rm Eoc}\limits}
\def\eoc{\mathop{\rm eoc}\limits}
\def\Pom{\mathop{\rm Pom}\limits}
\def\pom{\mathop{\rm pom}\limits}

\auteurcourant={DOMINIQUE FOATA AND GUO-NIU HAN}
\titrecourant={TREE CALCULUS FOR DIFFERENCE EQUATIONS}

\vglue.5cm

\rightline{December 21, 2012}
\bigskip\bigskip
\centerline{\bf Tree Calculus for
Bivariable Difference Equations}
\bigskip
\centerline{\sl Dominique Foata and Guo-Niu Han}
\footnote{}{
{\it Key words and phrases.} Tree Calculus, 
partial difference equations, strictly ordered binary trees,
end of minimal chain, parent maximum leaf, bivariate distributions,
Poupard triangle, tangent numbers.\par
{\it Mathematics Subject Classifications.} 
05A15, 05A30, 11B68, 33B10.}

\bigskip
{\narrower\narrower
\eightpoint
\noindent
{\bf Abstract}.\quad
Following Poupard's study of strictly ordered binary trees with respect to two parameters, namely, ``end of minimal chain'' and ``parent of maximum leaf'' a true Tree Calculus is being developed 
to solve a partial difference equation system and then make a joint study of those two statistics.  Their joint distribution is shown to  be symmetric and to be expressed in the form of an explicit three-variable generating function.

}
\bigskip\medskip
\centerline{\bf 1. Introduction}

\medskip
The triangle of numbers
{\eightpoint

$$
f=\matrice{&&&&f_0(1)\cr
&&&f_1(1)&f_1(2)&f_1(3)\cr
&&f_2(1)&f_2(2)&f_2(3)&f_2(4)&f_2(5)\cr
&f_{3}(1)&f_3(2)&f_3(3)&f_3(4)&f_3(5)&f_{3}(6)&f_{3}(7)\cr
f_{4}(1)&f_{4}(2)&f_4(3)&f_4(4)&f_4(5)&f_4(6)&f_{4}(7)
&f_{4}(8)&f_{4}(9)\cr
}\!=\!
\matrice{&&&&1\cr
&&&0&1&0\cr
&&0&1&2&1&0\cr
&0&4&8&10&8&4&0\cr
0&34&68&94&104&94&68&34&0\cr
}$$

}
\centerline{Table 1.1}

\smallskip
\noindent
appears in Sloane's On-Line Encyclopedia of Integer Sequences [Sl06] under reference A008301 and is called Poupard's triangle.
As shown by Christiane Poupard [Po89], $f=(f_{n}(m))$ 
$(n\ge 0,\,1\le m\le 2n+1)$ is the unique solution of the finite-difference equation system
$$
\Delta^2 f_n(m)+2\,f_{n-1}(m)
=0\quad(n\ge 1,\,1\le m\le 2n-1),\leqno(1.1)
$$
where $\Delta$ stands for the classical finite-difference
operator (see, e.g., [Jo39]) 
$$\leqalignno{
\Delta f_n(m)&:= f_n(m+1)-f_n(m),&(1.2)\cr
\noalign{\hbox{so that}} 
\Delta^2
f_n(m)&=f_n(m+2)-2f_n(m+1)+f_n(m),&(1.3)\cr}
$$
when taking $f_{0}(1)=1$; $f_{n}(1)=0$, 
\smash{$f_{n}(2)=\sum\limits_{m}f_{n-1}(m)$} 
for $n\ge 1$ as {\it initial values}. 
Let
$$
\leqalignno{ \noalign{\vskip-5pt}
\quad\tan u&=\sum_{n\ge 1} {u^{2n-1}\over
(2n-1)!}T_{2n-1}&(1.4)\cr
\noalign{\vskip-5pt}
&={u\over 1!}1+{u^3\over 3!}2+{u^5\over 5!}16+{u^7\over
7!}272+ {u^9\over 9!}7936+\cdots\cr}
$$

\goodbreak\noindent
be the Taylor expansion of $\tan u$, the coefficients $T_{2n+1}$
$(n\ge 0)$ being called the {\it tangent numbers} (see, e.g.,
[Ni23, p.~177-178], [Co74, p.~258-259]); Poupard further shows that each row sum $f_{n}(\brullet):=f_{n}(1)+f_{n}(2)+\cdots+f_{n}(2n+1)$ is equal to the {\it integer} $T_{2n+1}/2^{n}$ $(n\ge 0)$, that is, reporting to Table~1.1: 1, 1, 4, 34, 496,\dots\ 

\goodbreak
Finally, on the set ${\goth T}_{2n+1}$ of {\it strictly ordered binary trees} with $(2n+1)$ vertices (see Definition~1.3), she defines two statistics ``{\bf eoc}'' (``{\bf e}nd {\bf o}f minimal {\bf c}hain'') and ``{\bf pom}'' (``{\bf p}arent {\bf o}f the {\bf m}aximum leaf''), to show that both statistics ``eoc'' and ``pom+1'' are equally distributed on each set~${\goth T}_{2n+1}$, and furthermore,
$$
\#\{t\in {\goth T}_{2n+1}\!:\! \eoc (t)=k+1\}
\!=\!\#\{t\in {\goth T}_{2n+1}\!: \!\pom (t)=k\}=f_{n}(k)\leqno(1.5)
$$
for all $k$; in particular,
$\#{\goth T}_{2n+1}=T_{2n+1}/2^{n}$.

\medskip
The purpose of this paper is to calculate the {\it joint} distribution of the pair $(\eoc,\pom)$ on each set ${\goth T}_{2n+1}$ and to derive its properties, in particular its symmetry. To achieve this, we first introduce a sequence $(M_{n}=(f_{n}(m,k))$ of $(2n)\times(2n)$-matrices $(n\ge 1)$ with nonnegative integral entries, called a {\it Delta sequence}, defined by a system of finite difference equations, verifying certain initial conditions. Then, we show that each entry $f_{n}(m,k)$ is equal to the number of trees~$t$ from ${\goth T}_{2n+1}$ such that $\eoc(t)=m$ and $\pom(t)=k$.

\medskip
It is convenient to consider the following four triangles of each square $\{(m,k):1\le m,k\le 2n\}$:

$L_{n}^{(1)}:=\{2\le k+1\le m\le 2n-2\}$;\quad
$L_{n}^{(2)}:=\{4\le k+3\le m\le 2n\}$;

$U_{n}^{(1)}:=\{2\le m+1\le k\le 2n-2\}$;\quad
$U_{n}^{(2)}:=\{4\le m+3\le k\le 2n\}$.

\smallskip
\noindent
By convention, $f_{n}(m,k)\!:=\!0$ if $(m,k)\!\not\in\! [1,2n]\times[1,2n]$. The {\it partial differ\-ence operators} \smash{$\Deltaa_{m}$, $\Deltaa_{k}$,} act as follows on the entries of the matrices~$M_{n}$:
$$
\leqalignno{\Deltaa_{m}f_{n}(m,k)&:=f_{n}(m+1,k)-f_{n}(m,k);&(1.6)\cr
\Deltaa_{k}f_{n}(m,k)&:=f_{n}(m,k+1)-f_{n}(m,k).&(1.7)\cr}
$$
They serve to define the recurrence relations:
$$\leqalignno{
{\Deltaa_{m}}^2f_{n}(m,k)+2\,f_{n-1}(m,k)&=0\quad
((m,k)\in L_{n}^{(1)});&(R\,1)\cr
{\Deltaa_{k}}^2f_{n}(m,k)+2\,f_{n-1}(m,k)&=0\quad
((m,k)\in U_{n}^{(1)}).&(R\,2)\cr
}
$$
Finally, denote the {\it row} and {\it column sums} of~$M_{n}=(f_{n}(m,k))$ by
$$\leqalignno{
f_{n}(m,\brullet)&:=\sum\limits_{1\le k\le 2n}f_{n}(m,k);\quad
(1\le m\le 2n);\cr
f_{n}(\brullet,k)&:=\sum\limits_{1\le m\le 2n}f_{n}(m,k)\quad
(1\le k\le 2n).\cr}
$$ 

{\it Definition} 1.1.\quad
A sequence of matrices $(M_{n})$ $(n\ge 1)$, where each matrix $M_{n}=(f_{n}(m,k))\ (1\le k,m\le 2n)$ has nonnegative integral entries, having only 0's along its diagonal, and such that $M_{1}:=\petitematrice{0\;&0\cr 1\;&0\cr}$, is said to be a {\it Delta Sequence}, if for $n\ge 2$ both recurrence relations
$(R\,1)$ and $(R\,2)$ hold, together with the {\it initial conditions}:

\goodbreak
$(I\,1)$ for $n\ge 2$ the $(2n)$-th column, ${\rm Col}_{2n}$, of $M_{n}$ is
the zero-column; its the $(2n-1)$-st column, ${\rm Col}_{2n-1}$, is equal to
$$
f_{n-1}(1,\brullet),\ f_{n-1}(2,\brullet),\ 
\ldots,\ f_{n-1}(2n-2, \brullet),\ 0,\ 0,
$$
when read from top to bottom.

\smallskip
$(I\,2)$ the $(2n)$-th row of $M_{n}$ is also equal to
$$
f_{n-1}(1,\brullet),\ f_{n-1}(2,\brullet),\ 
\ldots,\ f_{n-1}(2n-2, \brullet),\ 0,\ 0,
$$
when read from left to right; its
$(2n-1)$-st row is equal to:
$$\displaylines{\quad
f_{n-1}(1,\brullet)+f_{n-1}(\brullet,1),\ 
f_{n-1}(2,\brullet)+f_{n-1}(\brullet,2),\ \ldots\ ,\hfill\cr
\hfill{}
f_{n-1}(2n-3,\brullet)+f_{n-1}(\brullet,2n-3),\ 
f_{n-1}(2n-2,\brullet)+f_{n-1}(\brullet,2n-2),\ 0,\ 0,\cr}$$
when read from left to right.

\medskip
In the above definition the entries of the matrix $M_{n}$ are derived from~$M_{n-1}$ by first applying rules $(I\,1)$ and $(I\,2)$ and letting the diagonal be null; then, starting from $m=1$ up to $m=2n-3$, for each $k$ from $2n-3$ down to $m+1$, evaluate $f_{n}(m,k)$ with equation $(R\,2)$: $f_{n}(m,k)-2f_{n}(m,k+1)+f_{n}(m,k+2)+2f_{n-1}(m,k)=0$, the coefficients $f_{n}(m,k+1)$, $f_{n}(m,k+2)$ and $f_{n-1}(m,k)$ being already calculated. Exchanging the roles of $m$ and~$k$ the upper entries are obtained by using equation $(R\,1)$. Accordingly, $(R\,1)$, $(R\,2)$, $(I\,1)$, $(I\,2)$ uniquely determine the Delta Sequence $(M_{n}$ $(n\ge 1)$.

\medskip
{\it Calculation of the first matrices}.\quad First,
$f_{1}(1,\brullet)=0$, $f_{1}(2,\brullet)=1$, so that
$$
M_{2}=\pmatrix{0&?&0&0\cr
?&0&1&0\cr
1&1&0&0\cr
0&1&0&0\cr}
$$
by rules $(I\,1)$ and $(I\,2)$.
The remaining entries are obtained by rule $(R\,2)$:
$f_{2}(1,2)-2f_{2}(1,3)+f_{2}(1,4)+2f_{1}(1,2)=f_{2}(1,2)-2\times 0+0+2\times 0=0$, so that $f_{2}(1,2)=0$; then, by rule~$(R\,1)$:
$f_{2}(2,1)-2f_{2}(3,1)+f_{2}(4,1)+2f_{1}(2,1)=f_{2}(2,1)-2\times 1+0+2\times 1=0$, so that $f_{2}(2,1)=0$. Thus,
$$
M_{2}=\pmatrix{0&0&0&0\cr
0&0&1&0\cr
1&1&0&0\cr
0&1&0&0\cr}
$$
and $f_{2}(1,\brullet)=0$, $f_{2}(2,\brullet)=1$, $f_{2}(3,\brullet)=2$, $f_{2}(4,\brullet)=1$. The next matrices are displayed in Fig.~1.2.

\goodbreak
\topinsert

$$M_{3}\!=\!\pmatrix{0&0&0&0&0&0\cr
0&0&1&2&1&0\cr
1&1&0&4&2&0\cr
2&3&4&0&1&0\cr
1&3&3&1&0&0\cr
0&1&2&1&0&0\cr}\  M_{4}\!=\!\pmatrix{0&0&0&0&0&0&0&0\cr
0&0&4&8&10&8&4&0\cr
4&4&0&16&20&16&8&0\cr
8&12&16&0&28&20&10&0\cr
10&18&24&28&0&16&8&0\cr
8&18&24&24&16&0&4&0\cr
4&12&18&18&12&4&0&0\cr
0&4&8&10&8&4&0&0\cr}\!\!.
$$

$$M_{5}=\pmatrix{
 0 &     0 &     0 &     0 &     0 &     0 &     0 &     0 &     0 &
  0 &   \cr
 0 &     0 &     34 &    68 &    94 &    104 &   94 &    68 &    34 &
  0 &   \cr
 34 &    34 &    0 &     136 &   188 &   208 &   188 &   136 &   68 &
  0 &   \cr
 68 &    102 &   136 &   0 &     274 &   296 &   262 &   188 &   94 &
  0 &   \cr
 94 &    162 &   222 &   274 &   0 &     352 &   296 &   208 &   104 &
  0 &   \cr
 104 &   198 &   276 &   330 &   352 &   0 &     274 &   188 &   94 &
  0 &   \cr
 94 &    198 &   282 &   330 &   330 &   274 &   0 &     136 &   68 &
  0 &   \cr
 68 &    162 &   240 &   282 &   276 &   222 &   136 &   0 &     34 &
  0 &   \cr
 34 &    102 &   162 &   198 &   198 &   162 &   102 &   34 &    0 &
  0 &   \cr
 0 &     34 &    68 &    94 &    104 &   94 &    68 &    34 &    0 &
  0 &  \cr
}
$$
\centerline{Fig. 1.2: the first matrices $M_{n}$}

\endinsert

\medskip
The previous definition of a Delta Sequence, based on the two relations $(R\,1)$, $(R\,2)$ and the two initial condtions $(I\,1)$, $(I\,2)$, can be symbolized by the square on the left in Fig.~1.3, as relation $(R\,1)$ (resp. $(R\,2)$) acts on the entries of the lower (resp. upper) entries of the matrix~$M_{n}$, and initial conditions $(I\,1)$ and $(I\,2)$ refer to the last two columns ${\rm Col}_{2n-1}$, ${\rm Col}_{2n}$ and rows ${\rm Row}_{2n-1}$, ${\rm Row}_{2n}$ of~$M_{n}$, respectively.
Other initial conditions will be stated in Section~8. We just mention a second one, materialized by the square on the right in Fig.~1.3.

$$
\vbox{\vskip1.6cm\offinterlineskip
\segment(0,0)\dir(1,0)\long{16}
\segment(0,2)\dir(1,0)\long{16}
\segment(0,0)\dir(0,1)\long{16}
\segment(16,0)\dir(0,1)\long{16}
\segment(0,16)\dir(1,0)\long{16}
\segment(14,0)\dir(0,1)\long{16}
\segment(0,16)\dir(1,-1)\long{16}
\centerput(15,8){$\longleftrightarrow$}
\centerput(8,-1){$\downarrow$}
\centerput(8,1.5){$\uparrow$}
\centerput(-4,-1){$\scriptstyle {\rm Row}_{2n}$}
\centerput(-6,2){$\scriptstyle {\rm Row}_{2n-1}$}
\centerput(10,17.5){$\scriptstyle {\rm Col}_{2n-1}$}
\centerput(20,17){$\scriptstyle {\rm Col}_{2n}$}
\centerput(22,8){$\scriptstyle (R\,2)$}
\centerput(8,-4){$\scriptstyle (R\,1)$}
}\hskip4cm\vbox{\vskip2cm\offinterlineskip
\segment(0,0)\dir(1,0)\long{16}
\segment(0,0)\dir(0,1)\long{16}
\segment(16,0)\dir(0,1)\long{16}
\segment(0,14)\dir(1,0)\long{16}
\segment(0,16)\dir(1,0)\long{16}
\segment(2,0)\dir(0,1)\long{16}
\segment(0,16)\dir(1,-1)\long{16}
\centerput(1,8){$\longleftrightarrow$}
\centerput(9,13){$\downarrow$}
\centerput(9,15.5){$\uparrow$}
\centerput(-7,8){$\scriptstyle (R\,4)$}
\centerput(9,20){$\scriptstyle (R\,3)$}
\centerput(-4,15){$\scriptstyle {\rm Row}_{1}$}
\centerput(-4,13){$\scriptstyle {\rm Row}_{2}$}
\centerput(-2,17){$\scriptstyle {\rm Col}_{1}$}
\centerput(4,17){$\scriptstyle {\rm Col}_{2}$}
}\hskip2cm
$$

\medskip
\centerline{Fig. 1.3: Definitions 1.1 and 1.2}

\medskip
{\it Definition} 1.2.\quad
A sequence of matrices $(M_{n})$ $(n\ge 1)$, where each matrix $M_{n}=(f_{n}(m,k))\ (1\le k,m\le 2n)$ has nonnegative integral entries, having only 0's along its diagonal, and such that $M_{1}:=\petitematrice{0\;&0\cr 1\;&0\cr}$, is said to be a {\it Gamma Sequence}, if for $n\ge 2$ both recurrence relations
$$
\leqalignno{\noalign{\vskip-5pt}
\qquad
{\Deltaa_{m}}^2f_{n}(m,k)+2\,f_{n-1}(m,k-2)&=0\quad
((k,m)\in U_{n}^{(2)});&(R\,3)\cr
\qquad
{\Deltaa_{k}}^2f_{n}(m,k)+2\,f_{n-1}(m-2,k)&=0\quad
((k,m)\in L_{n}^{(2)});&(R\,4)\cr
\noalign{\vskip-5pt}
}
$$
hold, together with the {\it initial conditions}:

\smallskip
$(I\,3)$ for $n\ge 2$ the first row, ${\rm Row}_{1}$, is the zero-row;
the second row, ${\rm Row}_{2}$, is equal to

\smallskip
$0, f_{n-1}(1,\brullet)(=0), f_{n-1}(2,\brullet),\ldots,
f_{n-1}(2n-2,\brullet),f_{n-1}(2n-1,\brullet)(=0)$;

\smallskip
\noindent
when read from left to right.

\smallskip
$(I\,4)$ the first column, ${\rm Col}_{1}$, of $M_{n}$ is also equal to

\smallskip
$0, f_{n-1}(1,\brullet)(=0), f_{n-1}(2,\brullet),\ldots,
f_{n-1}(2n-2,\brullet),f_{n-1}(2n-1,\brullet)(=0)$

\smallskip
\noindent
when read from top to bottom; the second column, ${\rm Col}_{2}$, is equal to
$$\displaylines{\noalign{\vskip-8pt}
\quad
0, 0, f_{n-1}(2,\brullet)+f_{n-1}(1,\brullet),\ 
f_{n-1}(3,\brullet)+f_{n-1}(2,\brullet),\ \ldots\hfill\cr
\hfill{}
f_{n-1}(2n-2,\brullet)+f_{n-1}(2n-3,\brullet),\ 
f_{n-1}(2n-2,\brullet),\cr
\noalign{\vskip-8pt}}$$
when read from left to right.

\medskip
Using the same reasoning as for Definition~1.1 it is seen that the Gamma Sequence is uniquely defined. The fact that Delta and Gamma Sequences are identical will be a consequence of the further theorems ({\it cf.} Section~6). Next, comes the combinatorial set-up on which all calculations will be made.

\medskip
{\it Definition} 1.3.\quad
An {\it $n$-labeled strictly ordered binary tree}
is defined by the following axioms:

(1) it is a {\it labeled} tree with $n$ nodes, labeled
$1,2,\ldots, n$; the node labeled~1 is called the {\it root};

\goodbreak
(2) each node has no child (it is then called a {\it leaf\/}), or two children, their order being immaterial
(it is then called an {\it interior node\/});

(3) when getting along each path from the root to each node,
the node labels are in {\it increasing} order.

\medskip
Each strictly ordered binary tree has an {\it odd} number of vertices, say, $2n+1$, with~$n$ interior nodes and $n+1$ leaves. Let ${\goth T}_{2n+1}$ denote the set of all strictly ordered binary trees with $(2n+1)$ nodes.
When giving an orientation (left or right) to each child of each of the~$n$ interior nodes in a tree~$t\in {\goth T}_{2n+1}$, we generate $2^n$ {\it planar} strictly ordered binary trees (also called ``arbres binaires croissants complets" by Viennot [Vi88, chap.~3, p.~111]). It is known that the latter are equidistributed with the {\it alternating permutations} of order~$(2n+1)$, so that their number is equal to the tangent number $T_{2n+1}$, a result that goes back to D\'esir\'e Andr\'e [An1879, An1881]. Accordingly,
$$
\#{\goth T}_{2n+1}=T_{2n+1}/2^n.\leqno(1.8)
$$

Let $t\in {\goth T}_{2n+1}$ $(n\ge 1)$. If a node labeled~$a$ has two children labeled $b$ and~$c$, define $\min a:=\min\{b,c\}$; if it has one child~$b$, let $\min a:=b$. The {\it minimal chain} of~$t$ is defined to be the sequence
$a_1\rightarrow a_2\rightarrow a_3\rightarrow\cdots
\rightarrow a_{j-1}\rightarrow a_j$, with the following properties: (i) $a_1=1$ is the label of the root; (ii) for
$i=1,2,\ldots,j-1$ the $(i+1)$-st term~$a_{i+1}$ is the label of an interior node and $a_{i+1}=\min a_i$; (iii) $a_j$ is the node of a leaf. Define the ``end of the minimal chain'' of~$t$ to be $\eoc(t):=a_{j}$. If the leaf with the maximum label~$(2n+1)$ is
incident to a node labeled~$k$, define its ``parent of the maximum leaf'' to be $\pom(t):=k$.

\newbox\boxarbre
\medskip
\setbox\boxarbre=\vbox{\vskip
22mm\offinterlineskip 
\centerput(30,21){$7$}
\centerput(40,21){$8$}
\centerput(10,21){$5$}
\centerput(20,21){$9$}
\centerput(15,12){$4$}
\centerput(37,14){$3$}
\centerput(0,11){$6$}
\centerput(30,7){$2$}\centerput(40,10){$\mapsto$}
\centerput(5,2){$1$}
\segment(0,10)\dir(1,-1)\long{5}
\segment(5,5)\dir(5,1)\long{25}
\segment(30,10)\dir(1,1)\long{5}
\segment(30,10)\dir(-3,1)\long{15}
\segment(15,15)\dir(1,2)\long{3}
\segment(15,15)\dir(-1,2)\long{3}
\segment(35,15)\dir(-1,2)\long{3}
\segment(35,15)\dir(1,2)\long{3}
}

\newbox\boxarbreb
\medskip
\setbox\boxarbreb=\vbox{\vskip
22mm\offinterlineskip 
\centerput(30,21){$9$}
\centerput(40,21){$7$}
\centerput(10,21){$4$}
\centerput(20,21){$8$}
\centerput(15,12){$3$}
\centerput(37,14){$6$}
\centerput(0,11){$5$}
\centerput(30,7){$2$}
\centerput(5,2){$1$}
\segment(0,10)\dir(1,-1)\long{5}
\segment(5,5)\dir(5,1)\long{25}
\segment(30,10)\dir(1,1)\long{5}
\segment(30,10)\dir(-3,1)\long{15}
\segment(15,15)\dir(1,2)\long{3}
\segment(15,15)\dir(-1,2)\long{3}
\segment(35,15)\dir(-1,2)\long{3}
\segment(35,15)\dir(1,2)\long{3}
}

$$
\box\boxarbre\hskip20mm\hskip28mm\box\boxarbreb
\hskip3.2cm
$$

\vskip-12pt
\centerline{Fig. 1.4}

\medskip
The minimal chain of the tree~$t$ displayed in
Fig.~1.4 on the left is $1\rightarrow 2\rightarrow 3\rightarrow 7$, so that $\eoc(t)=7$ and the parent of its maximum leaf (equal to $2n+1=9$) is $\pom(t)=4$. Also, the parent of the maximum leaf in the tree~$t'$ on the right is equal to $\pom(t')=6$. To go from $t$ to~$t'$ replace the labels $1, 2, 3, 7$ of the minimal chain by the labels $2-1, 3-1,7-1,9 = 1,2,6,9$ and change each other label~$a$ by $a-1$. We then have $\eoc(t)=\pom(t')+1=7$. This illustrates the construction of the bijection $t\mapsto t'$ of ${\goth T}_{2n+1}$ onto itself with the property $\eoc(t)=\pom(t')+1$, which was described in [Ha12]. Thus, the statistics~``$\eoc$'' and~``$\pom+1$'' are equidistributed on each ${\goth T}_{2n+1}$, as already mentioned in (1.5), a property first proved by Poupard [Po89] by analytic methods. The main results of this paper are stated next.

\proclaim Theorem 1.1. Let $(M_{n}=(f_{n}(m,k))$ $(n\ge 1)$ be the Delta sequence, as introduced in Definition~$1.1$. Then,
for all $n\ge 1$ and $1\le m,k\le 2n$
$$\#\{t\in{\goth T}_{2n+1}:\eoc(t)=m,\,\pom(t)=k\}=f_{n}(m,k).
\leqno(1.9)
$$

A first consequence of Theorem 1.1 and (1.8) is the identity:
$$\sum_{m,k} f_{n}(m,k)=T_{2n+1}/2^n\quad(n\ge 0).\leqno(1.10)
$$

\proclaim Theorem 1.2. Let $\bigl(M_{n}=(f_{n}(m,k))\ (1\le m,k\le 2n)\bigr)$ $(n\ge 1)$ be the Delta sequence. Then, the matrices $M_{n}$ are symmetric with respect to their counter-diagonals:
$$
f_{n}(m,k)=f_{n}(2n+1-k,2n+1-m)
\quad(1\le k,m\le 2n).\leqno(1.11)
$$

The proofs of those two theorems and also of further properties will be based on the geometric properties of the strictly ordered binary trees. To this end, we adopt the following notation and convention: for each triple $(n,m,k)$ let ${\goth T}_{2n+1,m,k}$  (resp. ${\goth T}_{2n+1,m,\brullet}$, resp. ${\goth T}_{2n+1,\brullet,k}$) denote the subset of ${\goth T}_{2n+1}$ of all trees~$t$ such that $\eoc(t)=m$ and $\pom(t)=k$ (resp. $\eoc(t)=m$, resp. $\pom(t)=k$). By convention, designate those families {\it and their cardinalities} by the same symbol and also the matrix of the {\it integers} ${\goth T}_{2n+1,m,k}$ by ${\rm Mat}({\goth T}_{2n+1})$. Our plan of action will be to show that the sequence $({\rm Mat}({\goth T}_{2n+1}))$ $(n\ge 1)$ {\it is identical} to the Delta Sequence.

In Sections~2--5 it will be shown that, when replacing each $f_{n}(m,k)$ by ${\goth T}_{2n+1,m,k}$ the initial conditions $(I\,1)$ and $(I\,2)$, the two finite-difference equations systems $(R\,1)$, $(R\,3)$, the two finite-difference equations systems $(R\,2)$, $(R\,4)$ and the initial conditions $(I\,3)$ and $(I\,4)$) hold. This will complete the proofs of Theorems~1.1 and~1.2, as done in Section~6. Further properties of the matrices ${\rm Mat}({\goth T}_{2n+1})$ (and then matrices~$M_{n}$) will be given in Section~7. Finally, several other equivalent definitions of the Delta sequence will be mentioned in Section~8. In Section~9 we conclude the paper by calculating the generating functions for the $f_{n}(m,k)$ in the following forms.

\proclaim Theorem 1.3. The triple-exponential generating function for the lower triangles of the matrices~$M_{n}$ is given by
$$
\displaylines{(1.12)\quad
\sum_{2\le k+1\le m\le 2n}
f_{n}(m,k){x^{m-k-1}\over (m-k-1)!}
{y^{k-1}\over (k-1)!}
{z^{2n-m}\over (2n-m)!}\hfill\cr
\hfill{}={\cos(\sqrt 2\, x)+\cos(\sqrt 2\, y)\,\cos(\sqrt 2\, z)\over
2\,\cos^2\Bigl(\displaystyle{x+y+z\over \sqrt2}\Bigr)}.\quad\cr
}$$

\proclaim Theorem 1.4. The triple exponential generating function for the upper triangles of the matrices~$M_{n}$ is given by
$$
\displaylines{(1.13)\quad
\sum_{2\le m+1\le k\le 2n-1}
f_{n}(m,k){x^{2n-k}\over (2n-k)!}
{y^{k-m-1}\over (k-m-1)!}
{z^{m-1}\over (m-1)!}\hfill\cr
\hfill{}=\sin(\sqrt 2\,x)\,
\sin(\sqrt 2\,z)\,
{1\over 2\,\displaystyle\cos^2\Bigl({x+y+z\over\sqrt 2}\Bigr)}.\quad\cr
}$$

\bigskip

\def\arbrea#1#2#3#4{\mathop{\hskip14pt
\vbox{\vskip1cm\offinterlineskip
\segment(0,0)\dir(0,1)\long{3}
\segment(0,3)\dir(1,1)\long{4}
\segment(0,3)\dir(-1,1)\long{4}
\rightput(-4,4.5){\hbox{$\scriptstyle#1$}}
\leftput(3.3,6.9){\hbox{$\scriptstyle#2$}}
\leftput(1,2){\hbox{$\scriptstyle#3$}}
\centerput(-4,7.5){\hbox{$\scriptstyle#4$}}
}\hskip14pt}\nolimits}

\def\carre{\mathop{\hbox{\kern5pt \vbox{
\offinterlineskip\segment(-1.5,-.5)\dir(0,1)\long{3}
\segment(-1.5,-0.5)\dir(1,0)\long{3}
\segment(1.5,-0.5)\dir(0,1)\long{3}
\segment(1.5,2.5)\dir(-1,0)\long{3}}\kern5pt }}}

\def\carresec{\mathop{\hbox{\kern5pt \vbox{
\offinterlineskip\segment(-2,-1)\dir(0,1)\long{4}
\segment(-2,-1)\dir(1,0)\long{4}
\segment(2,-1)\dir(0,1)\long{4}
\segment(2,3)\dir(-1,0)\long{4}
\segment(2,3)\dir(-2,1)\long{4}}\kern5pt }}}

\def\rond{\mathop{\raise.5pt\hbox{$\bigcirc$}}}

\def\rondbullet{\mathop{\hbox{\kern5pt \vbox{
\offinterlineskip
\centerput(0,2){\hbox{$\scriptstyle\bullet$}}
\centerput(0,0){\hbox{$\bigcirc$}}
\centerput(0,4){$\scriptstyle m$}
}}}}

\def\rondbullettwo{\mathop{\hbox{\kern5pt \vbox{
\offinterlineskip
\centerput(0,2){\hbox{$\scriptstyle\bullet$}}
\centerput(0,0){\hbox{$\bigcirc$}}
\centerput(-5,2){$\scriptstyle m-2$}
}}}}

\def\triang{\mathop{\raise.5pt\hbox{$\bigtriangledown$}}}
\def\triangup{\mathop{\raise.5pt\hbox{$\bigtriangleup$}}}

\def\triangbullet{\mathop{\hbox{\kern5pt \vbox{
\offinterlineskip
\centerput(0,2){\hbox{$\triang$}}
\centerput(1.7,3){$\scriptstyle\bullet$}
\centerput(3,4.4){$\scriptstyle m$}
}}}}

\def\arbrebb#1#2#3#4#5#6#7#8{\mathop{\hskip14pt
\vbox{\vskip1.2cm\offinterlineskip
\segment(0,0)\dir(0,1)\long{3}
\segment(0,3)\dir(1,1)\long{4}
\segment(0,3)\dir(-1,1)\long{8}
\segment(-4.5,7)\dir(1,1)\long{4}
\leftput(-10,11.5){$\scriptstyle#1$}
\leftput(2.5,7.6){\hbox{$\scriptstyle#2$}}
\leftput(-1.3,10.8){\hbox{$\scriptstyle#3$}}
\rightput(-4,4.5){\hbox{$\scriptstyle#4$}}
\rightput(-8,9){\hbox{$\scriptstyle#5$}}
\leftput(0.2,11){\hbox{$\scriptstyle#6$}}
\leftput(0.5,1.5){\hbox{$\scriptstyle#7$}}
\leftput(4,5){\hbox{$\scriptstyle#8$}}
}\hskip14pt}\nolimits}

\centerline{\bf 2. The initial conditions
$(I\,1)$ and $(I\,2)$}

\medskip
In this section and the next ones we make the convention that whenever a leaf is deleted from a tree, the edge linking the leaf to the tree is also deleted.

For verifying that the matrices 
${\rm Mat}({\goth T}_{2n+1})$ have only zero in their diagonals, it suffices to show that ${\goth T}_{2n+1,m,m}=\emptyset$ (or is equal to~0 with our convention). This is true, because if $t\in {\goth T}_{2n+1}$ and $\eoc(t)=\pom(t)=m$, the node $(2n+1)$ has a parent equal to~$m$. Consequently, $m$ cannot be the end of a minimal chain. Hence, the previous subset is empty.

\proclaim Theorem 2.1. The initial conditions $(I\,1)$ and $(I\,2)$ hold for the matrices 
${\rm Mat}({\goth T}_{2n+1})$ when $f_{n}(m,\brullet)$ and
$f_{n}(\brullet,k)$ are replaced by ${\goth T}_{2n+1,m,\brullet}$ and ${\goth T}_{2n+1,\brullet,k}$, respectively.

\proof
$(I\,1)$ First, the $(2n)$-th column of the matrix 
${\rm Mat}({\goth T}_{2n+1})$ has zero entries only, as
$\pom\le 2n-1$. Next,
each  tree from ${\goth T}_{2n+1,m,2n-1}$ 
$(1\le m\le 2n-2)$ must contain the subtree
$$\segment(1,5)\dir(1,-1)\long{4}
\segment(5,1)\dir(1,1)\long{4}
\centerput(-1,5){$\scriptstyle 2n$}
\centerput(5,-1){$\scriptstyle 2n-1$}
\centerput(13,5){$\scriptstyle 2n+1$}$$
Hence, ${\goth T}_{2n+1,2n-1,2n-1}$ is empty, for $(2n-1)$, being an interior node, cannot be the end of the minimal chain. Also,
${\goth T}_{2n+1,2n,2n-1}$ is empty, for the sibling of $(2n-1)$ is neccessarily less than $(2n-1)$, so that the minimal chain cannot go through $(2n-1)$ and reach~$(2n)$. Furthermore, ${\goth T}_{2n+1,1,2n-1}={\goth T}_{2n-1,1,\brullet}=0$, as $\eoc\ge 2$.

In the remaining cases, that is, $2\le m\le 2n-2$,
removing the two leaves $(2n)$, $(2n+1)$ transforms each tree from ${\goth T}_{2n+1,m,2n-1}$ onto a tree from ${\goth T}_{2n-1,m,\brullet}$ in a bijective manner. Such a transformation may be illustrated by the diagram:

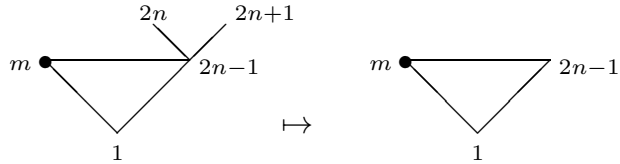
\def\arbrebbl{\mathop{\hskip14pt
\vbox{\vskip1.4cm\offinterlineskip
\segment(0,0)\dir(-1,1)\long{8}
\segment(0,0)\dir(1,1)\long{12}
\segment(8,8)\dir(-1,1)\long{4}
\segment(-8,8)\dir(1,0)\long{16}
\centerput(0,-3){$\scriptstyle1$}
\leftput(9,6.5){\hbox{$\scriptstyle 2n-1$}}
\rightput(-9.5,7){\hbox{$\scriptstyle m$}}
\centerput(-8,7){$\bullet$}
\centerput(4,12.5){$\scriptstyle 2n$}
\leftput(13,12.5){\hbox{$\scriptstyle 2n+1$}}
}\hskip14pt}\nolimits}

\def\arbrebbm{\mathop{\hskip14pt
\vbox{\vskip1.4cm\offinterlineskip
\segment(0,0)\dir(-1,1)\long{8}
\segment(0,0)\dir(1,1)\long{8}
\segment(-8,8)\dir(1,0)\long{16}
\centerput(0,-3){$\scriptstyle1$}
\leftput(9,6.5){\hbox{$\scriptstyle 2n-1$}}
\rightput(-9.5,7){\hbox{$\scriptstyle m$}}
\centerput(-8,7){$\bullet$}
}\hskip14pt}\nolimits}
$$\arbrebbl\kern 35pt\mapsto\kern35pt \arbrebbm$$

\goodbreak\smallskip\noindent
Hence, the $(2n-1)$-st column of the matrix ${\rm Mat}({\goth T}_{2n+1})$ reads:
$$
{\goth T}_{2n-1,1,\brullet},{\goth T}_{2n-1,2,\brullet},\ldots,
{\goth T}_{2n-1,2n-2,\brullet},\ 0,\ 0.\leqno(2.1)
$$
from top to bottom.

\def\arbrebbl{\mathop{\hskip14pt
\vbox{\vskip1.4cm\offinterlineskip
\segment(0,0)\dir(-1,1)\long{8}
\segment(0,0)\dir(1,1)\long{12}
\segment(8,8)\dir(-1,1)\long{4}
\segment(-8,8)\dir(1,0)\long{16}
\centerput(0,-3){$\scriptstyle1$}
\leftput(9,6.5){\hbox{$\scriptstyle k$}}
\rightput(4,11.5){\hbox{$\bullet$}}
\centerput(4,13.5){$\scriptstyle 2n$}
\leftput(13,12.5){\hbox{$\scriptstyle 2n+1$}}
}\hskip14pt}\nolimits}

\def\arbrebbm{\mathop{\hskip14pt
\vbox{\vskip1.4cm\offinterlineskip
\segment(0,0)\dir(-1,1)\long{8}
\segment(0,0)\dir(1,1)\long{8}
\segment(-8,8)\dir(1,0)\long{16}
\centerput(0,-3){$\scriptstyle1$}
\leftput(9,6.5){\hbox{$\scriptstyle k$}}
\rightput(8.5,7){\hbox{$\bullet$}}
}\hskip14pt}\nolimits}

$(I\,2)$\quad
For the $(2n)$-th row of the matrix ${\rm Mat}({\goth T}_{2n+1})$ note that ${\goth T}_{2n+1,2n,1}=0$ for $n\ge 2$. When $k\ge 2$ each tree from ${\goth T}_{2n+1,2n,k}$ must contain the subtree
$$\segment(1,5)\dir(1,-1)\long{4}
\segment(5,1)\dir(1,1)\long{4}
\centerput(-2,5){$\scriptstyle 2n$}
\centerput(1,4.5){$\bullet$}
\centerput(5,-1){$\scriptstyle k$}
\centerput(13,5){$\scriptstyle 2n+1$}$$
By $(I\,1)$ we then have: ${\goth T}_{2n+1,2n,2n-1}={\goth T}_{2n+1,2n,2n}=0$. For the remaining cases $2\le k\le 2n-2$ we can set up a bijection of ${\goth T}_{2n+1,2n,k}$ onto ${\goth T}_{2n-1,k,\brullet}$ by removing the two leaves $2n$ and $(2n+1)$, as illustrated by the next diagram.
$$\arbrebbl\kern 35pt\mapsto\kern35pt \arbrebbm$$

\smallskip
\noindent
Note that the node~$k$ becomes the end of the minimal chain. Thus, the $(2n)$-th row of the matrix ${\goth T}_{2n+1}$ is also equal to (2.1) read from left to right.

Finally, consider the $(2n-1)$-st row of ${\goth T}_{2n+1}$. In an obvious manner, ${\goth T}_{2n+1,2n-1,2n-1}={\goth T}_{2n+1,2n-1,2n}=0$ . When $1\le k\le 2n-2$, the trees from the {\it sets}
${\goth T}_{2n+1,2n-1,k}$ fall into two categories ${\goth T}_{2n+1,2n-1,k}^{I}$ and ${\goth T}_{2n+1,2n-1,k}^{II}$. In the first category the trees contain the subtree

\bigskip\medskip\smallskip\noindent
$
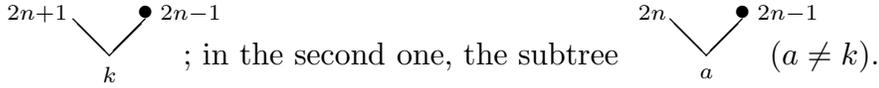
\qquad 
\segment(1,5)\dir(1,-1)\long{4}
\segment(5,1)\dir(1,1)\long{4}
\centerput(-3,5){$\scriptstyle 2n+1$}
\centerput(5,-2){$\scriptstyle k$}
\centerput(9,5){$\bullet$}
\centerput(14,5){$\scriptstyle 2n-1$}$
\kern1.2cm; in the second one, the subtree \quad
$\segment(1,5)\dir(1,-1)\long{4}
\segment(5,1)\dir(1,1)\long{4}
\centerput(-1,5){$\scriptstyle 2n$}
\centerput(5,-1.5){$\scriptstyle a$}
\centerput(14,5){$\scriptstyle 2n-1$}
\centerput(9,5){$\bullet$}\kern1.2cm(a\not=k)$.

\medskip
First, note that 
${\goth T}_{2n+1,2n-1,1}^{I}={\goth T}_{2n-1,1,\bullet}=0$. When $2\le k\le 2n-2$, removing the two leaves $(2n+1)$, $(2n-1)$ and replacing the node label $(2n)$ by $(2n-1)$ maps ${\goth T}_{2n+1,2n-1,k}^{I}$ onto ${\goth T}_{2n-1,k,\bullet}$ in a bijective manner. 

When $1\le k\le 2n-2$, removing the two leaves $(2n)$, $(2n-1)$, and replacing the node label $(2n+1)$ by $(2n-1)$ maps ${\goth T}_{2n+1,2n-1,k}^{II}$ onto ${\goth T}_{2n-1,\brullet,k}$ in a bijective manner. Thus, the $(2n-1)$-st row of ${\goth T}_{2n+1}$ reads
$$\displaylines{\quad
{\goth T}_{2n-1,1,\brullet}+{\goth T}_{2n-1,\brullet,1},\  
{\goth T}_{2n-1,2,\brullet}+{\goth T}_{2n-1,\brullet,2},\ 
\ldots\ ,\hfill\cr
\hfill{}{\goth T}_{2n-1,2n-3,\brullet}+{\goth T}_{2n-1,\brullet,2n-3},\ 
{\goth T}_{2n-1,2n-2,\brullet}+{\goth T}_{2n-1,\brullet,2n-2},\ 
0,\ 0.\quad\qed\cr}
$$

\medskip
\centerline{\bf 3. Tree Calculus for the relations $(R\,1)$ and
$(R\,3)$}

\medskip

In the following Tree Calculus subtrees (possibly leaves) are indicated by the symbols ``$\rond$," ``$\triang$", or ``$\carre$." The end of the minimal chain in each tree is represented by a bullet ``$\bullet$.'' Letters occurring below or next to subtrees are labels of their roots. For instance, the symbols
$$
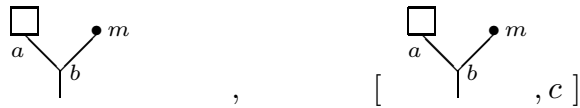
\arbrea{a}{\bullet\, m}b{\carre}\qquad\qquad,\qquad\qquad
[\quad\arbrea{a}{\bullet\, m}b{\carre}\quad,c\ ]
$$
designate the {\it families} of all trees~$t$ from the underlying set ${\goth T}_{2n+1}$ having a node labeled~$b$ [in short, a node~$b$], parent of both a
subtree of root~$a$ and the leaf~$m$, which is also the end of the minimal chain; moreover, the symbol on the right
has the further property that the node labeled~$c$ does not belong, 
{\it either} to the subtree of root~$b$, {\it or} to  the path going from root~1 to~$b$. In the sequel,
the letter ``$m$" is always used to designate the end of the minimal chain, unless explicitly indicated by a letter next to $\bullet$.

\def\arbrebb#1#2#3#4#5#6#7#8{\mathop{\hskip14pt
\vbox{\vskip1.2cm\offinterlineskip
\segment(0,0)\dir(0,1)\long{3}
\segment(0,3)\dir(1,1)\long{4}
\segment(0,3)\dir(-1,1)\long{8}
\segment(-4.5,7)\dir(1,1)\long{4}
\leftput(-10,11.5){$\scriptstyle#1$}
\leftput(2.5,7.6){\hbox{$\scriptstyle#2$}}
\leftput(-1.3,10.8){\hbox{$\scriptstyle#3$}}
\rightput(-4,4.5){\hbox{$\scriptstyle#4$}}
\rightput(-8,9){\hbox{$\scriptstyle#5$}}
\leftput(0.2,11){\hbox{$\scriptstyle#6$}}
\leftput(0.5,1.5){\hbox{$\scriptstyle#7$}}
\leftput(4,5){\hbox{$\scriptstyle#8$}}
}\hskip14pt}\nolimits}

\def\arbrebbj#1#2#3#4#5#6#7#8{\mathop{\hskip14pt
\vbox{\vskip1.2cm\offinterlineskip
\segment(0,0)\dir(0,1)\long{3}
\segment(0,3)\dir(-1,1)\long{4}
\segment(0,3)\dir(1,1)\long{8}
\segment(4.5,7)\dir(-1,1)\long{4}
\leftput(-10,11.5){$\scriptstyle#1$}
\leftput(7.5,9){\hbox{$\scriptstyle#2$}}
\leftput(-1.3,10.8){\hbox{$\scriptstyle#3$}}
\centerput(-4,6.5){\hbox{$\scriptstyle#4$}}
\rightput(10,11.5){\hbox{$\scriptstyle#5$}}
\leftput(-1.5,11.5){\hbox{$\scriptstyle#6$}}
\rightput(-4,4.5){\hbox{$\scriptstyle#7$}}
\leftput(4,5){\hbox{$\scriptstyle#8$}}
}\hskip14pt}\nolimits}

Our Tree Calculus consists of two steps: (a) decomposing the sets ${\goth T}_{2n+1,m,k}$ into smaller subsets by considering the mutual positions of the nodes $m$, $(m+1)$, $(m+2)$ (resp. $k$, $(k+1)$, $(k+2)$); (b) setting up bijections between those subsets by a simple display of certain subtrees, as done in (3.1).

\smallskip
For instance,
$$
C_{3}:=\qquad\arbrebb{\lower2pt\hbox{$\;\bullet$}}{\carre}%
{\rond}{m}{m+1}{}{}{m+2}\qquad
{\rm and}\qquad D_{1}:=
\qquad\arbrebbj{}{m+2}{}%
{\bullet}{\carre}{\rond}{m}{m+1}\qquad\leqno(3.1)
$$
may be regarded as two subsets of ${\goth T}_{2n+1}$.
To each pair
$(\hbox{$\vrule height8pt width0pt\carre \atop m+2$},\raise2pt\hbox{$\rond$})$, the nodes of~``$\rond$''  being all greater than or equal to $(m+2)$, there correspond a unique tree from $C_{3}$  and a unique tree from~$D_{1}$. This clearly defines a bijection of $C_{3}$ onto~$D_{1}$.

\smallskip

Those two principles (a) and (b) will be applied in the proofs of the next two theorems~3.1 and~4.1.

\proclaim Theorem 3.1. If $(m,k)$ belongs to
$L_{n}^{(1)}\cup U_{n}^{(2)}
=\{2\le k+1\le m\le 2n-2\}\cup
\{4\le m+3\le k\le 2n\}$, then
$$\leqalignno{\noalign{\vskip-10pt}
{\Deltaa_{m}}^2\,
{\goth T}_{2n+1,m,k}+2\quad\arbrea{m+1}{}{m}{\lower2pt\hbox{$\bullet$}}
\centerput(1,5){$\scriptstyle m+2$}&=0,&(3.2)\cr}$$
with the understanding that the second term on the left-hand side represents twice the set of all trees from ${\goth T}_{2n+1,m+1,k}$
with the further property that~$m$ is the parent of both~$(m+1)$ and $(m+2)$.

\proof
The decomposition
$$\leqalignno{
{\goth T}_{2n+1,m,k}
&=\quad\arbrea{m}{\lower4.5pt\hbox{\kern-3pt$\rond\atop
m+1$}}{} {\lower2pt\hbox{$\bullet$}}\quad
+\quad\arbrea{m}{\rond}{}{\lower2pt\hbox{$\bullet$}}
\cr 
}
$$
means that in each tree from
${\goth T}_{2n+1,m,k}$  the node $(m+1)$ is, or is not, the
sibling of the leaf~$m$. In the next decomposition the node~$m$
is, or is not, the parent of the leaf~$(m+1)$:
$$
{\goth T}_{2n+1,m+1,k}
=\quad\arbrea{m+1}{\rond}{m} {\lower2pt\hbox{$\bullet$}}\quad
+\quad\arbrea{m+1}{\rond}{b}{\lower2pt\hbox{$\bullet$}}\qquad(b\not=m).
$$
Under the transposition $(m,m+1)$ the node labeled~$k$ remains unaffected, because $k\le m-1$ if $(m,k)\in L_{n}^{(1)}$ and $m+3\le k$ if $(m,k)\in U_{n}^{(2)}$, so that the parent of $(2n+1)$ remains~$k$. Thus, the transposition establishes a one-to-one correspondence between the two second terms. Hence,
$$
\leqalignno{
\Deltaa_m{\kern-3pt}^2&\,
{\goth T}_{2n+1,m,k}\cr
&=
({\goth T}_{2n+1,m+2,k}-{\goth T}_{2n+1,m+1,k})
-({\goth T}_{2n+1,m+1,k}-{\goth T}_{2n+1,m,k})\cr
\noalign{\medskip}
&=\quad\arbrea{m+2}{\raise1.5pt\hbox{\kern-3pt$\rond$}}{m+1} {\lower1.5pt\hbox{$\bullet$}}\qquad
-\qquad\arbrea{m+1}{\lower4.5pt\hbox{\kern-5pt$\rond\atop
\ m+2$}}{} {\lower2pt\hbox{$\bullet$}}\quad
-\qquad\arbrea{m+1}{\lower0pt\hbox{$\rond
$}}{m}{\lower2pt\hbox{$\bullet$}}\qquad
+\qquad\arbrea{m}{\lower4.5pt\hbox{\kern-5pt$\rond\atop
\ m+1$}}{} {\lower2pt\hbox{$\bullet$}}\quad\cr 
&:=\kern18pt A\kern30pt  -\kern30pt B
\kern28pt -\kern29pt C\kern20pt +\kern18pt D\quad
.\cr}$$
Depending on the mutual positions of nodes~$m$, $(m+1)$ and $(m+2)$ the further decompositions prevail, as again, 
$k$ remaining still attached to $(2n+1)$:
$$
\leqalignno{
A&=\qquad\arbrebb{\lower2pt\hbox{$\;\bullet$}}{\carre}%
{\rond}{m+1}{m+2}{}{m}{}\quad
+\quad[\quad\arbrebb{\lower2pt\hbox{$\;\bullet$}}{\carre}%
{\rond}{m+1}{m+2}{}{}{}\quad,m]:=A_{1}+A_{2};\cr
\noalign{\medskip}
B&=\quad \arbrea{m+1}{}{m}{\lower2pt
\hbox{$\bullet$}}\centerput(1,5){$\scriptstyle m+2$}\quad
+\quad
\arbrebbj{}{}{}{\textstyle\bullet}{\rond}{\carre}{m+1}{m+2}
\centerput(-1.5,2){$\scriptstyle m$}\quad
+
\quad[\quad\arbrea{m+1}{\lower4.5pt\hbox{\kern-3pt$\rond\atop
m+2$}}{} {\lower2pt\hbox{$\bullet$}}\quad,m]
:=B_{1}+B_{2}+B_{3};\cr
\noalign{\medskip}
C&=\quad \arbrea{m+1}{}{m}{\lower2pt
\hbox{$\bullet$}}\centerput(1,5){$\scriptstyle m+2$}\quad
+\quad
\arbrebbj{}{}{}{\textstyle\bullet}{\rond}{\carre}{m+1}{m+2}
\centerput(-1.5,2){$\scriptstyle m$}\quad
+\qquad\arbrebb{\lower2pt\hbox{$\;\bullet$}}{\carre}%
{\rond}{m}{m+1}{}{}{m+2}\quad
+\quad[\quad\arbrebb{\lower2pt\hbox{$\;\bullet$}}{\carre}%
{\rond}{m}{m+1}{}{}{}\quad,m+2]\cr
&:=C_{1}+C_{2}+C_{3}+C_{4};\cr
\noalign{\medskip}
D&=\quad\arbrebbj{}{m+2}{}%
{\bullet}{\carre}{\rond}{m}{m+1}\qquad
+\quad[\quad\arbrea{m}{\kern-5pt \lower3pt\hbox{$\rond\atop \ m+1$}}{}{\lower2pt
\hbox{$\bullet$}}\quad,m+2]:=D_{1}+D_{2}.\cr}
$$

The permutation ${\  m\quad m+1\ m+2\choose m+2\ m\quad m+1}$
establishes a one-to-one correspondence between~$A_{2}$ and~$C_{4}$ (resp. between~$B_{3}$ and $D_{2}$). On the other hand,
$A_{1}=B_{2}+C_{2}=2\,B_{2}$, because the two subtrees in $B_{2}=C_{2}$ are incident to the same node, while in~$A_{1}$ they are incident to two different nodes. Finally, $C_{3}=D_{1}$ as explained in (3.1). Making all the proper cancellations in the sum
$\Deltaa_m{\kern-3pt}^2\,{\goth T}_{2n+1,m,k}
=A-B-C+D$ we get identity (3.2).\qed

\proclaim Corollary 3.2. The relations $(R\,1)$ and $(R\,3)$ hold for the matrices ${\rm Mat}({\goth T}_{2n+1})$, that is,
$$\leqalignno{
\Deltaa_m{\kern-3pt}^2\,{\goth T}_{2n+1,m,k}
+2\,{\goth T}_{2n-1,m,k}&=0,\quad{\sl if}\ (m,k)\in L_{n}^{(1)};&(3.3)\cr
\Deltaa_m{\kern-3pt}^2\,{\goth T}_{2n+1,m,k}
+2\,{\goth T}_{2n-1,m,k-2}&=0,\quad{\sl if}\ (m,k)\in U_{n}^{(2)}.&(3.4)\cr}
$$

\proof For (3.3) change the second term in (3.2) as follows: remove the two leaves $(m+1)$, $(m+2)$, and subtract~2 from all the remaining nodes greater than~$(m+2)$: the term  becomes ${\goth T}_{2n-1,m,k}$, as the node $(2n+1)$ becomes $(2n-1)$ and is still linked to~$k$. For (3.4) do the same changes, but this time, as $m+3\le k$, the edge going from~$k$ to $(2n+1)$ becomes and edge going from $(k-2)$ to $(2n-1)$.\qed

\bigskip
\centerline{\bf 4. Tree Calculus for the relations $(R\,2)$ and
$(R\,4)$}

\medskip
\proclaim Theorem 4.1. If $(m,k)$ belongs to $L_{n}^{(2)}
\cup U_{n}^{(1)}=\{4\le k+3\le m\le 2n\}\cup
\{2\le m+1\le k\le 2n-2\}$, then
$$\leqalignno{\noalign{\bigskip}
\Deltaa_k{\kern-3pt}^2\,{\goth T}_{2n+1,m,k}
+2\,[\quad\arbrebb{}{\carre}%
{}{k+1}{k+2}{2n+1}{k}{}\quad,m]
+2\,
\quad\arbrebb{\rondbullet}%
{k+2}{}{k+1}{}{2n+1}{k}{}
&=0,&(4.1)\cr
}$$
with the understanding that the second term on the left-hand side is twice the number of all trees from ${\goth T}_{2n+1,m,k+1}$ with the further property that $(k+2)$ is a leaf incident to~$(k+1)$, itself incident to~$k$, the end~$m$ of the minimum chain
being outside the subtree of root~$k$.

\goodbreak
\noindent{\it Proof}.
First,
$$\leqalignno{\noalign{\vskip-10pt}
{\goth T}_{2n+1,m,k}
&=[\arbrea{k+1}{\ 2n+1}k{\rond},m]+
[\arbrea{}{\ 2n+1}k{\rond}, m,k+1]
+\arbrea{k+1}{\ 2n+1}k{\rondbullet\ }
+\quad[\arbrea{}{\ 2n+1}k{\rondbullet\ },k+1]
\cr
&:=A_{1}+A_{2}+A_{3}+A_{4},\cr
}
$$
meaning that each tree from ${\goth T}_{2n+1,m,k}$ has one of
the four forms: either $k+1$ is incident to~$k$, or not, and $m$ is outside or not the subtree of root~$k$; furthermore, the leaf~$m$ is the end of the minimal chain. 

\goodbreak

Using the same dichotomy,
$$\leqalignno{\noalign{\vskip5pt}
{\goth T}_{2n+1,m,k+1}
&=[\quad\arbrebb{\rond}{\carre}{}{k+1}{}{2n+1}{k}{},m]
+[\arbrea{}{\ 2n+1}{k+1}{\rond}\quad,m, k]
+\quad\arbrebb{\rondbullet}{\carre}{}{k+1}{}{2n+1}{k}{}+\quad[\arbrea{}{\ 2n+1}{k+1}{\rondbullet\ },\ k]
\cr
\noalign{\medskip}
&:=B_{1}+B_{2}+B_{3}+B_{4}.\cr
\noalign{\medskip\smallskip}}$$
Consider the subsets 
$A'_{4}:=\quad\arbrebb{\rondbullet}%
{\lower5pt\hbox{$\,\carre\atop  k+1$}}{}{k}{}{2n+1}{}{}$ of $A_{4}$ and
$B'_{2}:=\quad\arbrebb{\rond}%
{\lower5pt\hbox{$\carre\atop k$}}{}{k+1}{}{2n+1}{}{\raise12pt\hbox{$\scriptstyle\ \bullet\ m$}}$ of $B_{2}$.
The transposition $(k,k+1)$ maps $A_{4}\setminus A'_{4}$ onto $B_{4}$ and $A_{2}$ onto $B_{2}\setminus B'_{2}$ in a bijective manner.

 Hence,
$$\leqalignno{\noalign{\vskip-5pt}
{\goth T}_{2n+1,m,k+1}-{\goth T}_{2n+1,m,k}
&=(B_{1}-A_{1})+((B_{2}-B'_{2}-A_{2})+B'_{2})\cr
&\quad{}+(B_{3}-A_{3})+(B_{4}-(A_{4}-A'_{4})-A'_{4})\cr
&=B_{1}-A_{1}+B'_{2}+B_{3}-A_{3}-A'_{4}\cr
\noalign{\medskip}
&\kern-2cm=[\quad\arbrebb{\rond}{\carre}{}{k+1}{}{2n+1}{k}{},m]
-[\quad\arbrea{k+1}{\ 2n+1}k{\rond},m]\ 
+\quad\arbrebb{\rond}%
{\lower5pt\hbox{$\carre\atop k$}}{}{k+1}{}{2n+1}{}{\raise12pt\hbox{$\scriptstyle\ \bullet\ m$}}
\cr
\noalign{\bigskip}
&\kern.5cm{}
+\quad\arbrebb{\rondbullet}{\carre}{}{k+1}{}{2n+1}{k}{}-\quad\arbrea{k+1}{\ 2n+1}k{\rondbullet\ }
-\qquad\arbrebb{\rondbullet}%
{\lower5pt\hbox{$\,\carre\atop  k+1$}}{}{k}{}{2n+1}{}{}\quad
;\cr}
$$

\goodbreak
\noindent and also
$$\leqalignno{
{\goth T}_{2n+1,m,k+2}-{\goth T}_{2n+1,m,k+1}
&=[\quad\arbrebb{\rond}{\carre}{}{k+2}{}{2n+1}{k+1}{}\ ,m]
-[\quad\arbrea{k+2}{\ 2n+1}{k+1}{\rond}\ ,m]
+\quad\arbrebb{\rond}%
{\lower5pt\hbox{$\carre\atop k+1$}}{}{k+2}{}{2n+1}{}{\raise12pt\hbox{$\scriptstyle\kern6pt \bullet\ m$}}
\cr
\noalign{\bigskip}
&\kern1cm{}+\quad\arbrebb{\rondbullet}{\carre}{}{k+2}{}{2n+1}{k+1}{}\ -\quad\arbrea{k+2}{\  2n+1}{k+1}{\rondbullet\ }\quad
-\qquad\arbrebb{\rondbullet}%
{\lower5pt\hbox{$\,\carre\atop  k+2$}}{}{k+1}{}{2n+1}{}{}
\quad\cr
&:=D_{1}-C_{1}+D'_{2}+D_{3}-C_{3}-C'_{4}.\cr}
$$
Thus, 
$$\leqalignno{
\Deltaa_k{\kern-3pt}^2\,{\goth T}_{2n+1,m,k}&=\bigl({\goth T}_{2n+1,m,k+2}-{\goth T}_{2n+1,m,k+1}\bigr)
-\bigl({\goth T}_{2n+1,m,k+1}-{\goth T}_{2n+1,m,k}\bigr)\cr
&=D_{1}-C_{1}+D'_{2}+D_{3}-C_{3}-C'_{4}\cr
&\qquad{}
-B_{1}+A_{1}-B'_{2}-B_{3}+A_{3}+A'_{4}.\cr
}
$$

\def\arbrebbc#1#2#3#4#5#6#7#8#9{\mathop{\hskip14pt
\vbox{\vskip1.6cm\offinterlineskip
\segment(4,-1)\dir(0,1)\long{3}
\segment(4,2)\dir(1,1)\long{4}
\segment(0,6)\dir(1,1)\long{4}
\segment(0,6)\dir(-1,1)\long{8}
\segment(0,6)\dir(1,-1)\long{4}
\segment(-4.5,10)\dir(1,1)\long{4}
\leftput(-10,14.5){$\scriptstyle#1$}
\leftput(2.5,10.6){\hbox{$\scriptstyle#2$}}
\leftput(-1.3,13.8){\hbox{$\scriptstyle#3$}}
\rightput(-4,7.5){\hbox{$\scriptstyle#4$}}
\rightput(-8,12){\hbox{$\scriptstyle#5$}}
\leftput(0.2,14){\hbox{$\scriptstyle#6$}}
\rightput(-0.5,4.5){\hbox{$\scriptstyle#7$}}
\leftput(4,8){\hbox{$\scriptstyle#8$}}
\centerput(8,6.5){$\scriptstyle#9$}
}\hskip14pt}\nolimits}

\def\arbrebbz#1#2#3#4#5#6#7#8#9{\mathop{\hskip14pt
\vbox{\vskip2cm\offinterlineskip
\segment(0,6)\dir(1,1)\long{10}
\segment(6,12)\dir(-1,1)\long{4}
\segment(0,6)\dir(-1,1)\long{8}
\segment(0,6)\dir(0,-1)\long{4}
\segment(-4.5,10)\dir(1,1)\long{4}
\leftput(-10,14.5){$\scriptstyle#1$}
\centerput(6,9){\hbox{$\scriptstyle#2$}}
\leftput(-1.3,13.5){\hbox{$\scriptstyle#3$}}
\rightput(-4,7.5){\hbox{$\scriptstyle#4$}}
\rightput(1,15){\hbox{$\scriptstyle#5$}}
\leftput(1,16.5){\hbox{$\scriptstyle#6$}}
\rightput(2,3.5){\hbox{$\scriptstyle#7$}}
\leftput(3.5,9.5){\hbox{$\scriptstyle#8$}}
\centerput(10,16.5){$\scriptstyle#9$}
}\hskip14pt}\nolimits}

\vskip-8pt
The further decompositions of the components of the previous sum depend on the mutual positions of the nodes~$k$, $(k+1)$, $(k+2)$. First, evaluate the subsum: $S_{1}:=D_{1}-C_{1}-B_{1}+A_{1}$ using the decompositions:
$$
\leqalignno{
D_{1}=[\qquad\arbrebb{\rond}{\carre}{}{k+2}{}{2n+1}{k+1}{}\quad
,\ m]
&=[\qquad\arbrebb{\rond}{\carre}{}{k+2}{}{2n+1}{k+1}{}\quad
,\ m,k\ ]+
[\quad\arbrebbc{\rond}{\carre}{}{k+2}{}{2n+1}{k+1}{}{\triang}
\raise3pt\hbox{$\scriptstyle k$}\quad,m]\cr
&:=D_{1,1}+D_{1,2};\cr
\noalign{\medskip}
C_{1}=[\quad\arbrea{k+2}{\,2n+1}{k+1}{\rond}\quad
,\ m]
&=[\quad\arbrea{k+2}{\,2n+1}{k+1}{\rond}\quad
,\ m,k\ ]\quad
+\quad[\quad\arbrebb{\rond}{\lower2.5pt\hbox{$\scriptstyle\
\carre$}}{}{k+1}{k+2}{2n+1}{k}{}\quad,\ m]\cr
&:=C_{1,1}+C_{1,2};\cr
\noalign{\medskip}
B_{1}=[\quad\arbrebb{\rond}{\carre}{}{k+1}{}{2n+1}{k}{},\ m]
&=[\quad\arbrebb{\rond}{\carre}{}{k+1}{}{2n+1}{k}{},\ m,k+2]
\quad
+[\quad\arbrebb{\lower2pt\hbox{$\;$}}{\carre}%
{}{k+1}{k+2}{2n+1}{k}{},m]\quad
\cr
&\kern-2cm{}
+\quad
[\arbrebbc{\rond}{\lower2pt\hbox{$\;$}}{\raise2pt
\hbox{$\kern-2pt
\triang$}}{k+2}{}{}{k+1}{\raise5pt\hbox{$\scriptstyle\kern4pt
2n+1$}}{\carre}\raise3pt\hbox{$\scriptstyle k$}
{\raise10pt\hbox{$\scriptstyle\kern-8pt $}}\quad,m]
+\quad[\quad \arbrebb{\carre}{\lower3pt\hbox{$\ $}}{}{k+1}{}{2n+1}{k}{k+2}\quad,m]
+\quad[\quad
\lower9pt\hbox{$\arbrebbz{\carre}{\quad k+2}{}{k+1}{2n+1}{\rond}{k}%
{}{\triang}$}\quad,m]\cr
&:=B_{1,1}+B_{1,2}\cr
&\kern1cm{}
+B_{1,3}+B_{1,4}+B_{1,5};\cr
\noalign{\medskip}
A_{1}=[\quad\arbrea{k+1}{\,2n+1}{k}{\rond}
,\ m]
&=[\quad\arbrea{k+1}{\,2n+1}{k}{\rond}
,\ m,k+2\ ]+[\qquad\arbrebb{\rond}{\lower2.5pt\hbox{$\scriptstyle\
$}}{\carre}{k+1}{k+2}{}{k}{2n+1}\quad,m]\cr
&:=A_{1,1}+A_{1,2}.\cr
}
$$
Also, let
\vskip-15pt
$$D'_{1,1}:=
\arbrebbc{\rond}{\carre}{}{k+2}{}{2n+1}{k+1}{}%
{\kern-3pt\lower5pt\hbox{${\triangbullet}$}}
{\kern8pt\raise12pt\hbox{$\scriptstyle k$}}\quad;
\qquad C'_{1,1}:=\quad\arbrebb{\rond}%
{\lower6pt\hbox{$\triangbullet$}}{}{k+1}{k+2}%
{\raise2pt\hbox{$\scriptstyle 2n+1$}}%
{{\kern8pt\raise9pt\hbox{$\scriptstyle k$}}}{}\quad.$$

The permutation ${\  k\quad k+1\ k+2\choose k+2\ k\quad k+1}$
maps~$D_{1,1}\setminus D'_{1,1}$ onto~$B_{1,1}$ and~$C_{1,1}
\setminus C'_{1,1}$ onto~$A_{1,1}$. 
Hence, $D_{1,1}=B_{1,1}+D'_{1,1}$,\ $C_{1,1}=A_{1,1}+C'_{1,1}$.
Moreover, $D_{1,2}=2\,B_{1,3}$,\ $C_{1,2}=A_{1,2}$,\ $B_{1,2}=B_{1,4}$,\ $B_{1,3}=B_{1,5}$. Altogether,
$S_{1}=D_{1}-C_{1}-B_{1}+A_{1}=(B_{1,1}+D'_{1,1}+2\,B_{1,3})
-(A_{1,1}+C'_{1,1}+A_{1,2})
-(B_{1,1}+B_{1,2}+B_{1,3}+B_{1,2}+B_{1,3})
+(A_{1,1}+A_{1,2})$.
Thus,
$$
S_{1}=-2\,B_{1,2}+D'_{1,1}-C'_{1,1}.\leqno(4.2)
$$

Next, evaluate the sum
$S_{2}:=D'_{2}+D_{3}-C_{3}-C'_{4}
-B'_{2}-B_{3}+A_{3}+A'_{4}$ by decomposing each of its components. We have:
$$
\leqalignno{\noalign{\medskip}
D'_{2}=\quad\arbrebb{\rond}%
{\lower5pt\hbox{$\carre\atop k+1$}}{}{k+2}{}{2n+1}{}{\raise12pt\hbox{$\scriptstyle\kern6pt \bullet\ m$}}
&=\quad[\quad\arbrebb{\rond}%
{\lower5pt\hbox{$\carre\atop k+1$}}{}{k+2}{}{2n+1}{}{\raise12pt\hbox{$\scriptstyle\kern6pt \bullet\ m$}},\quad k]
+\quad\arbrebb{\rond}%
{\lower5pt\hbox{$\carre\atop k+1$}}{}{k+2}{}{2n+1}{k}{\raise12pt\hbox{$\scriptstyle\kern6pt \bullet\ m$}}\cr
&:=D'_{2,1}+D'_{2,2};\cr
\noalign{\medskip}
D_{3}=\quad\arbrebb{\rondbullet}{\carre}{}{k+2}{}{2n+1}{k+1}{}
&=\quad[\quad\arbrebb{\rondbullet}{\carre}{}{k+2}{}{2n+1}{k+1}{},\quad k]+\quad\arbrebbc{\rondbullet}{\carre}{}{k+2}{}{2n+1}{k+1}{}%
{\triang}{\scriptstyle k}\qquad\cr 
&:=D_{3,1}+D_{3,2};\cr
}$$

\goodbreak
$$\leqalignno{
C_{3}=\quad\arbrea{k+2}{\  2n+1}{k+1}{\rondbullet\ }\quad
&=\quad[\quad\arbrea{k+2}{\  2n+1}{k+1}{\rondbullet\ }\quad,k]
+\quad\arbrebb{\rondbullet}{\carre}{}{k+1}{k+2}{2n+1}{k}{}\quad\cr
&:=C_{3,1}+C_{3,2};\cr
C'_{4}=\qquad\arbrebb{\rondbullet}%
{\lower5pt\hbox{$\,\carre\atop  k+2$}}{}{k+1}{}{2n+1}{}{}
&=\quad[\quad\arbrebb{\rondbullet}%
{\lower5pt\hbox{$\,\carre\atop  k+2$}}{}{k+1}{}{2n+1}{}{}
,k] \quad+\quad\arbrebb{\rondbullet}%
{k+2}{}{k+1}{}{2n+1}{k}{}
\quad+\quad
\lower9pt\hbox{$\arbrebbz{\rondbullet}{\quad k+2}{}{k+1}{2n+1}{\carre}{k}
{}{\lower1pt\hbox{$\triang$}\ }$}   \cr
&:=C'_{4,1}+C'_{4,2}+C'_{4,3};\cr
B'_{2}=
\quad\arbrebb{\rond}%
{\lower5pt\hbox{$\carre\atop k$}}{}{k+1}{}{2n+1}{}{\raise12pt\hbox{$\scriptstyle\ \bullet\ m$}}
&=\quad[\quad\arbrebb{\rond}%
{\lower5pt\hbox{$\carre\atop k$}}{}{k+1}{}{2n+1}{}{\raise12pt\hbox{$\scriptstyle\ \bullet\ m$}}\ ,k+2]
\ +\ 
\quad\arbrebb{\rond}%
{\lower5pt\hbox{$\carre\atop k$}}{}{k+1}{k+2}{2n+1}{}{\raise12pt\hbox{$\scriptstyle\ \bullet\ m$}}\quad+\qquad
\lower9pt\hbox{$\arbrebbz{\carre}{\quad k}{}{k+1}{2n+1}{\rond}{}
{}{\lower9pt\hbox{$\triangbullet\atop \quad k+2$}\ }$}\cr
&:=B'_{2,1}+B'_{2,2}+B'_{2,3};\cr
\cr
B_{3}=
\quad\arbrebb{\rondbullet}{\carre}{}{k+1}{}{2n+1}{k}{}
&=\quad[\quad\arbrebb{\rondbullet}{\carre}{}{k+1}{}{2n+1}{k}{},\  k+2]
+\ \arbrebb{\rondbullet}{\carre}{}{k+1}{k+2}{2n+1}{k}{}
\ +\quad\arbrebb{\rondbullet}%
{k+2}{}{k+1}{}{2n+1}{k}{}
\quad+\quad
\lower9pt\hbox{$\arbrebbz{\rondbullet}{\quad k+2}{}{k+1}{2n+1}{\carre}{k}
{}{\lower1pt\hbox{$\triang$}\ }$} 
\cr
&:=B_{3,1}+B_{3,2}+B_{3,3}+B_{3,4};\cr
\noalign{\bigskip\smallskip}
A_{3}=\quad\arbrea{k+1}{\ 2n+1}k{\rondbullet\ }
&=\quad[\quad\arbrea{k+1}{\ 2n+1}k{\rondbullet\ },\quad k+2]
\ 
+\qquad\arbrebb{\rondbullet}{2n+1}{}{k+1}{k+2}{\kern-2pt\carre }{k}{}\cr
&:=A_{3,1}+A_{3,2};\cr
\noalign{\bigskip}
A'_{4}=\quad\arbrebb{\rondbullet}%
{\lower5pt\hbox{$\,\carre\atop  k+1$}}{}{k}{}{2n+1}{}{}\quad
&=\quad[\quad\arbrebb{\rondbullet}%
{\lower5pt\hbox{$\,\carre\atop  k+1$}}{}{k}{}{2n+1}{}{}\ ,k+2]
\ +\ \arbrebb{\rondbullet}%
{\lower5pt\hbox{$\,\carre\atop  k+1$}}{}{k}{k+2}{2n+1}{}{}\quad
+\quad
\lower9pt\hbox{$\arbrebbz{\rondbullet}{\quad k+1}{}{k}{2n+1}{\carre}{}
{}{\lower7pt\hbox{$\triang\atop \quad k+2$}\ }$}\cr
&:=A'_{4,1}+A'_{4,2}+A'_{4,3}.\cr
}$$

\noindent
Within the sum $S_{2}$ there are numerous cancellations we now describe.

(a) {\it Components of the form $[t,k]$ or $[t,k+2]$, where $t$ is a subtree, whose root is labeled}.\quad There are four of them:
$D_{3,1}$, $-C_{3,1}$, $-B_{3,1}$, $A_{3,1}$.
Consider the subsets:
$$
B_{3,1,1}:=\quad\arbrebbc{\rondbullet}{\carre}{}{k+1}{}{2n+1}{k}{}{\lower 4pt\hbox{$\triang\atop \quad\  k+2$}}{}\quad;\quad
A_{3,1,1}:=
\quad\arbrebb{\rondbullet}{}{}{k}{k+1}{2n+1}{\raise 12pt\hbox{$\triang\atop \kern6pt  k+2$}}{}
\quad;
$$
of $B_{3,1}$ and $A_{3,1}$, respectively. The permutation ${\  k\quad k+1\ k+2\choose k+2\ k\quad k+1}$
maps $D_{3,1}$ onto $B_{3,1}\setminus B_{3,1,1}$ and
$C_{3,1}$ onto $A_{3,1}\setminus A_{3,1,1}$. 
Hence, $D_{3,1}-C_{3,1}-B_{3,1}+A_{3,1}
=(B_{3,1}-B_{3,1,1})-(A_{3,1}-A_{3,1,1})-B_{3,1}+A_{3,1}=
-B_{3,1,1}+A_{3,1,1}$.

\smallskip
(b) {\it Components of the form $[t,k]$ or $[t,k+2]$, where $t$ is a subtree, whose root is {\it not} labeled}.\quad There are four of them:
$D'_{2,1}$, $-C'_{4,1}$, $-B'_{2,1}$, $A'_{4,1}$. Again,
the permutation ${\  k\quad k+1\ k+2\choose k+2\ k\quad k+1}$
maps $D'_{2,1}$ onto $B'_{2,1}$, and~$C'_{4,1}$ onto $A'_{4,1}$. Hence,
$D'_{2,1}-B'_{2,1}=-C'_{4,1}+A'_{4,1}=0$. Their sum vanish.

\smallskip
(c) {\it Components represented by a tree $t$, whose root is {\it unlabeled}}.\quad There are four of them: $-B'_{2,2}$, $-B'_{2,3}$, $-A'_{4,2}$, $A'_{4,3}$. As $B'_{2,2}=A'_{4,2}$, the contribution of those components to~$S_{2}$ is then $-B'_{2,3}+A'_{4,3}$.

\smallskip
(d) {\it Components represented by a tree $t$, whose root is {\it labeled}}.\quad There are nine of them:
$D'_{2,2}$, $D_{3,2}$, $-C_{3,2}$, $-C'_{4,2}$, $-C'_{4,3}$
 $-B_{3,2}$, $-B_{3,3}$, $-B_{3,4}$,
$A_{3,2}$.
By simply comparing the subtree contents we have:
$D'_{2,2}-C_{3,2}=
-B_{3,2}+A_{3,2}=0$, $D_{3,2}-(C'_{4,3}+B_{3,4})=0$ and
$C'_{4,2}=B_{3,3}$. The contribution of those terms is then
$-2\,C'_{4,2}$. 

\smallskip
Hence, $S_{1}+S_{2}=(-2\,B_{1,2}+D'_{1,1}-C'_{1,1})+\bigl((-B_{3,1,1}+A_{3,1,1})+(-B'_{2,3}+A'_{4,3})+(-2\,C'_{4,2})\bigr)$. As 
$D'_{1,1}=B'_{2,3}$, $C'_{1,1}=A_{3,1,1}$ and 
$B_{3,1,1}=A'_{4,3}$, we get $S_{1}+S_{2}=-2\,B_{1,2}-2\,C'_{4,2}$, that is,
$$
\Deltaa_k{\kern-3pt}^2\,{\goth T}_{2n+1,m,k}
-2\,[\quad\arbrebb{\lower2pt\hbox{$\;$}}{\carre}%
{}{k+1}{k+2}{2n+1}{k}{},m]\quad-2\,
\quad\arbrebb{\rondbullet}%
{k+2}{}{k+1}{}{2n+1}{k}{}
.\qed
$$

\proclaim Corollary 4.2. The relations $(R\,2)$ and $(R\,4)$ hold for the matrices ${\rm Mat}({\goth T}_{2n+1})$, that is,
$$\leqalignno{
\Deltaa_k{\kern-3pt}^2\,{\goth T}_{2n+1,m,k}
+2\,{\goth T}_{2n-1,m-2,k}&=0,\quad{\sl if}\ (m,k)\in L_{n}^{(2)};&(4.3)\cr
\Deltaa_k{\kern-3pt}^2\,{\goth T}_{2n+1,m,k}
+2\,{\goth T}_{2n-1,m,k}&=0,\quad{\sl if}\ (m,k)\in U_{n}^{(1)}.&(4.4)\cr}
$$

\proof
When $(m,k)$ belongs to $L_{n}^{(2)}$, the second term in 
(4.1) is in\medskip\noindent bijection with $2[\quad\arbrea{2n-1}{\carre}{k} {},m-2]$ and the third one with
$2 \arbrea
{\raise8pt\hbox{$\rondbullettwo$}\kern-0pt}
{\ 2n-1}{k} {\raise8pt\hbox{}}$, that is, the set of trees in which the end $m-2$ of the minimal chain is outside (resp. inside) the subtree of root~$k$. 
The sum of those two terms is then $2\,{\goth T}_{2n-1,m-2,k}$.

\goodbreak
When $(m,k)$ belongs to $U_{n}^{(1)}$, the third term of (4.1) vanishes and the second one is in bijection with 
$2[\quad\arbrea{2n-1}{\carre}{k} {},m]$, equal to
$2\,{\goth T}_{2n-1,m,k}$.\qed

\bigskip\smallskip
\centerline{\bf 5. The initial conditions
$(I\,3)$ and $(I\,4)$}

\medskip

\def\arbreag{\mathop{\hskip14pt
\vbox{\vskip1cm\offinterlineskip
\segment(0,0)\dir(1,1)\long{12}
\segment(0,8)\dir(1,0)\long{8}
\segment(0,0)\dir(-1,1)\long{4}
\segment(4,4)\dir(-1,1)\long{4}
\segment(8,8)\dir(-1,1)\long{4}
\centerput(-4,3.5){$\bullet$}
\centerput(0,-3){$\scriptstyle 1$}
\centerput(8,5){$\scriptstyle k$}
\centerput(4,1){$\scriptstyle 3$}
\centerput(-6,4){$\scriptstyle 2$}
\centerput(15,9.5){$\scriptstyle 2n+1$}
\centerput(4,12.5){$\carre$}
}}}

\def\arbreah{\mathop{\hskip14pt
\vbox{\vskip1cm\offinterlineskip
\segment(0,0)\dir(1,1)\long{8}
\segment(-4,4)\dir(1,0)\long{8}
\segment(0,0)\dir(-1,1)\long{4}
\segment(4,4)\dir(-1,1)\long{4}
\centerput(0,-3){$\scriptstyle 1$}
\centerput(11.5,6){$\scriptstyle 2n-1$}
\centerput(7,2.5){$\scriptstyle k-2$}
\centerput(0,8.5){$\carre$}
}}}

\def\arbreai{\mathop{\hskip14pt
\vbox{\vskip1cm\offinterlineskip
\segment(0,0)\dir(1,1)\long{8}
\segment(0,8)\dir(1,0)\long{8}
\segment(0,0)\dir(-1,1)\long{4}
\segment(4,4)\dir(-1,1)\long{4}
\centerput(8,7.5){$\bullet$}
\centerput(0,-3){$\scriptstyle 1$}
\centerput(9,6){$\scriptstyle m$}
\centerput(4,1){$\scriptstyle 2$}
\centerput(-9,4){$\scriptstyle 2n+1$}
}}}

\def\arbreaj{\mathop{\hskip14pt
\vbox{\vskip1cm\offinterlineskip
\segment(0,0)\dir(1,1)\long{4}
\segment(-4,4)\dir(1,0)\long{8}
\segment(0,0)\dir(-1,1)\long{4}
\centerput(0,-3){$\scriptstyle 1$}
\centerput(8,2.5){$\scriptstyle m-1$}
\centerput(4,3){$\bullet$}
}}}

\proclaim Property 5.1. Initial conditions 
$(I\,3)$ and $(I\,4)$ hold for the matrices ${\rm Mat}({\goth T}_{2n+1})$.

\proof $(I\,3)$\quad
The first row of each matrix ${\rm Mat}({\goth T}_{2n+1})$ is obviously the zero-row, as~1 can never be the end of the minimal chain. For the second row note that for $n\ge 2$ the set ${\goth T}_{2n+1,2,k}$ is empty when $k=2$ or $2n$. Also the set ${\goth T}_{2n+1,2,1}$ is empty, for 2 and $(2n+1)$ can be both children of the root only when $n=1$.

Let $3\le k\le 2n-1$. As illustrated by the diagram

\medskip\smallskip
$$
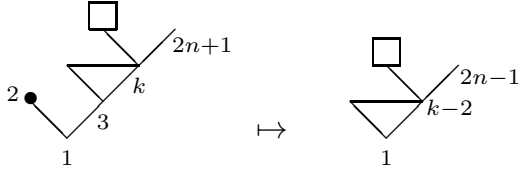
\arbreag\kern2cm\mapsto\kern.5cm \arbreah\kern1cm$$

\medskip
\noindent
each tree from ${\goth T}_{2n+1,2,k}$
is transformed into a tree from ${\goth T}_{2n-1,\brullet,k-2}$
by deleting the two nodes~2 and~1 and reducing the remaining nodes by~2. This transformation is obviously a bijection. Thus, ${\goth T}_{2n+1,2,k} =
{\goth T}_{2n-1,\brullet,k-2}$ for $3\le k\le 2n-1$. The second row of~${\rm Mat}({\goth T}_{2n+1})$ is then equal to the sequence 
$$\displaylines{\rlap{(5.1)}\hfill
0,0,{\goth T}_{2n-1,\brullet,1}, 
{\goth T}_{2n-1,\brullet,2},\ldots,{\goth T}_{2n-1,\brullet,2n-3},
{\goth T}_{2n-1,\brullet,2n-2}(=0);\hfill\cr
\noalign{\hbox{which is also equal to}}
\rlap{(5.2)}\hfill
0, {\goth T}_{2n-1,1,\brullet}(=0),
{\goth T}_{2n-1,2,\brullet}, \ldots,
{\goth T}_{2n-1,2n-2,\brullet},0,\hfill\cr}
$$
by Poupard's result (1.5) (also by our combinatorial proof sketched in the Introduction).

$(I\,4)$\quad The set ${\goth T}_{2n+1,m,1}$ is empty when $m=1,2$ or~$2n$. When $3\le m\le 2n-1$ the diagram

$$\arbreai\kern1.5cm\mapsto\kern.5cm  \arbreaj\kern1cm$$

\medskip
\noindent
serves to illustrate the transformation that maps
each tree from ${\goth T}_{2n+1,m,1}$ onto a tree
from ${\goth T}_{2n-1,m-1,\brullet}$, by deleting the two nodes $(2n\!+\!1)$ and~1, and reducing the remaining nodes by~1. Thus, the first column of~$M_{n}$ is equal to sequence (5.2),
when read from top to bottom.

For the second column we first note that
${\goth T}_{2n+1,m,2}$ is empty when $m=1,2$ and $n\ge 2$. When $m=3$, the mapping
$$\mathop{\hskip14pt
\vbox{\vskip1cm\offinterlineskip
\segment(0,0)\dir(1,1)\long{8}
\segment(0,0)\dir(-1,1)\long{4}
\segment(4,4)\dir(-1,1)\long{4}
\centerput(8,7.5){$\bullet$}
\centerput(0,-3){$\scriptstyle 1$}
\centerput(4,1){$\scriptstyle 2$}
\centerput(8,5){$\scriptstyle 3$}
\centerput(0,9){$\scriptstyle 2n+1$}
\centerput(-4,4.5){$\rond$}
}}\kern1cm\mapsto\kern.5cm
\mathop{\hskip14pt
\vbox{\vskip1cm\offinterlineskip
\segment(0,0)\dir(1,1)\long{4}
\segment(0,0)\dir(-1,1)\long{4}
\centerput(0,-3){$\scriptstyle 1$}
\centerput(8,3){$\scriptstyle 2n+1$}
\centerput(-4,4.5){$\rond$}
}}\kern1cm
$$
\smallskip\noindent
shows that ${\goth T}_{2n+1,3,2}$ is in bijection with
${\goth T}_{2n-1,\brullet,1}={\goth T}_{2n-1,2,\brullet}$.
When $m=4$, the following decomposition prevails
$$\leqalignno{
{\goth T}_{2n+1,4,2}&=
\mathop{\hskip14pt
\vbox{\vskip1.5cm\offinterlineskip
\segment(0,0)\dir(1,1)\long{12}
\segment(0,0)\dir(-1,1)\long{4}
\segment(4,4)\dir(-1,1)\long{4}
\segment(8,8)\dir(-1,1)\long{4}
\centerput(12,11){$\bullet$}
\centerput(0,-3){$\scriptstyle 1$}
\centerput(14,10.5){$\scriptstyle 4$}
\centerput(4,1){$\scriptstyle 2$}
\centerput(8,5){$\scriptstyle 3$}
\centerput(0,9){$\scriptstyle 2n+1$}
\centerput(-4,4.5){$\rond$}
\centerput(4,12.5){$\carre$}
}}\kern1.5cm+\kern.5cm
\mathop{\hskip14pt
\vbox{\vskip1cm\offinterlineskip
\segment(0,0)\dir(1,1)\long{8}
\segment(0,0)\dir(-1,1)\long{4}
\segment(4,4)\dir(-1,1)\long{4}
\centerput(8,7){$\bullet$}
\centerput(0,-3){$\scriptstyle 1$}
\centerput(4,1){$\scriptstyle 2$}
\centerput(8,5){$\scriptstyle 4$}
\centerput(0,9){$\scriptstyle 2n+1$}
\centerput(-4,4.5){$\rond$}
\centerput(-4,1){$\scriptstyle 3$}
}}\cr
\noalign{\medskip}
\noalign{\hbox{the two sets on the right-hand side
being in bijection, respectively, with}}
&=
\mathop{\hskip14pt
\vbox{\vskip1.5cm\offinterlineskip
\segment(0,0)\dir(1,1)\long{8}
\segment(0,0)\dir(-1,1)\long{4}
\segment(4,4)\dir(-1,1)\long{4}
\centerput(8,7){$\bullet$}
\centerput(0,-3){$\scriptstyle 1$}
\centerput(4,1){$\scriptstyle 2$}
\centerput(8,5){$\scriptstyle 3$}
\centerput(0,8.5){$\carre$}
\centerput(-4,4.5){$\rond$}
}}\kern1cm+\kern.5cm
\mathop{\hskip14pt
\vbox{\vskip1cm\offinterlineskip
\segment(0,0)\dir(1,1)\long{4}
\segment(0,0)\dir(-1,1)\long{4}
\centerput(0,-3){$\scriptstyle 1$}
\centerput(6,3){$\scriptstyle 2$}
\centerput(4,3.5){$\bullet$}
\centerput(-4,4.5){$\rond$}
\centerput(-4,1){$\scriptstyle 3$}
}}\qquad,\cr
\noalign{\hbox{that is,}}
&={\goth T}_{2n-1,3,\brullet}
+{\goth T}_{2n-1,2,\brullet}.\cr
}
$$

To prove the identity 
${\goth T}_{2n+1,m,2}={\goth T}_{2n-1,m-1,\brullet}
+{\goth T}_{2n-1,m-2,\brullet}
={\goth T}_{2n+1,m,1}
+{\goth T}_{2n+1,m-1,1}$ (by (5.2) and $(I\,4)$)
for $5\le m\le 2n$  proceed by induction on~$m$ using relation $(R\,1)$ already proved in (3.3). Write
$$
\eqalignno{
{\goth T}_{2n+1,m,2}&=2\,{\goth T}_{2n+1,m-1,2}\!-\!{\goth T}_{2n+1,m-2,2}\!-\!2\,{\goth T}_{2n-1,m-2,2}\qquad\qquad&\hbox{[by (3.3)]}\cr
2\,{\goth T}_{2n+1,m-1,2}
&=2\,{\goth T}_{2n+1,m-1,1}
+2\,{\goth T}_{2n+1,m-2,1}&\hbox{[by induction on $m$]}\cr
{}-{\goth T}_{2n+1,m-2,2}
&=-{\goth T}_{2n+1,m-2,1}
-{\goth T}_{2n+1,m-3,1}&\hbox{[by induction on $m$]}\cr
{}-2\,{\goth T}_{2n-1,m-2,2}
&=-2\,{\goth T}_{2n-1,m-2,1}
-2\,{\goth T}_{2n-1,m-3,1}&\cr
&&\hbox{[by induction on $n$ and $m$]}\cr
{}2\,{\goth T}_{2n+1,m-1,1}&-{\goth T}_{2n+1,m-2,1}-2\,{\goth T}_{2n-1,m-2,1}=
{\goth T}_{2n+1,m,1}&\hbox{[by (3.3)]}\cr
{}2\,{\goth T}_{2n+1,m-2,1}&-{\goth T}_{2n+1,m-3,1}-2\,{\goth T}_{2n-1,m-3,1}=
{\goth T}_{2n+1,m-1,1}.&\hbox{[by (3.3)]}\cr
}
$$
By summing those six equations we get
${\goth T}_{2n+1,m,2}={\goth T}_{2n+1,m,1}+{\goth T}_{2n+1,m-1,1}$, also equal to
${\goth T}_{2n-1,m-1,\brullet}
+{\goth T}_{2n-1,m-2,\brullet}$ by (5.2).\qed

\medskip
{\it Remark}.\quad It would be interesting to make up a proof of Property~5.1 that would have no recourse to a recurrence argument for $m\ge 5$ as above.

\vfill\eject
\centerline{\bf 6. Proofs of Theorems 1.1 and 1.2}

\medskip
Taking Property 2.1, Corollaries 3.2 and 4.2, Property~5.1
into account we conclude that the sequence of matrices
${\rm Mat}({\goth T}_{2n+1})$ is both a Delta and Gamma sequence. Those sequences are then identical and we may write
$$
f_{n}(m,k)={\goth T}_{2n+1,m,k}\leqno(6.1)
$$
for all $n,m,k$. This completes the proof of Theorem~1.1.

\medskip
We now exploit
the properties of the strictly ordered binary trees to prove that the matrices $M_{n}$ are symmetric with respect to their counter-diagonals (Theorem~1.2). 
First, the symmetry property is banal for $M_{1}$, $M_{2}$. For $n\ge 3$ consider the NE- and SW-corners
$$\leqalignno{
\petitematrice{f_{n}(1,2n-1)&f_{n}(1,2n)\cr
f_{n}(2,2n-1)&f_{n}(2,n)\cr}&=
\petitematrice{0&0\cr
f_{n}(2,2n-1)&0\cr}\cr
\noalign{\medskip}
\petitematrice{f_{n}(2n-1,1)&f_{n}(2n-1,2)\cr
f_{n}(2n,1)&f_{n}(2n,2)\cr}
&=\petitematrice{f_{n}(2n-1,1)&f_{n}(2n-1,2)\cr
0&f_{n}(2n,2)\cr}\cr}
$$
of the matrix $M_{n}$. As $f_{n}(2n-1,1)=f_{n}(2n,2)=
{\goth T}_{2n-3}=T_{2n-3}/2^{n-2}$
(by combining (1.11), $(I\,2)$, $(I\,3)$, (5.1) and (5.2)), both corners are symmetric with respect to their counter-diagonals [in short, counter-symmetric].

Let us prove that the {\it upper part} of the matrix $M_{n}$ is counter-symmetric and for $i=1,2,\ldots, n-1$ adopt the notation:
$$\leqalignno{{\rm Row}_{i}
&=\{(i,i+1),(i,i+2),\ldots,(i,2n-i),(i,2n-i+1)\};\cr
{\rm Col}_{2n-i+1}
&=\{(i,2n-i+1),(i+1,2n-i+1),\cr
&\kern1cm{}\ldots,(2n-i-1,2n-i+1),(2n-i,2n-i+1)\}.\cr
}
$$ 
Note that ${\rm Row}_i$ and ${\rm Col}_{2n-i+1}$ have the cell 
$(i,2n-i+1)$, belonging to the counter-diagonal, in common.
There is nothing to prove for the pairs $(m,k)$ along the
counter-diagonal and also for the entries from ${\rm Row}_{1}$ and ${\rm Col}_{2n}$, which are all zero.

Let $(m_0,k_0)$ belong to ${\rm Row}_j$ for some~$j$ such
that $2\le j\le n-1$. Further, assume that (1.12) holds for all 
$(m,k)\in{\rm Row}_{1}\cup\cdots\cup{\rm Row}_{j-1}$
and all $(m,k)\in {\rm Row}_j$ lying on the right of
$(m_0,k_0)$, not including $(m_0,k_0)$, that is, $m=j$ and $k>k_0$. By symmetry, (1.12) also holds for all $(m,k)\in {\rm
Col}_{2n}\cup\cdots\cup{\rm Col}_{2n-j+2}$ and all $(m,k)\in
{\rm Col}_{2n-j+1}$ lying above
$(2n+1-k_0,2n+1-m_0)$ not including the latter pair.

Now, the following relations hold:
$$
\displaylines{\ 
f_{n}(m_0,k_0)\!=\!2f_{n}(m_0,k_0+1)-f_{n}(m_0,k_0+2)
-2\,f_{n-1}(m_0,k_0)\hfill\hbox{[by $(R\,2)$]}\cr
\  f_{n-1}(m_0,k_0)=f_{n-1}(2n-1-k_0,2n-1-m_0),
\hfill\hbox{[by induction on~$n$]}\cr
\ f_{n}(m_0,k_0+1)=f_{n}(2n-k_0,2n+1-m_0),\hfill\cr
\ f_{n}(m_0,k_0+2)=f_{n}(2n-1-k_0,2n+1-m_0).
\hfill\hbox{[both by the local induction]}\cr
}
$$
Therefore,
$$
\displaylines{
f_{n}(m_0,k_0)-2\,f_{n}(2n-k_0,2n+1-m_0)\hfill
\cr
\hfill{}+f_{n}(2n-1-k_0,2n+1-m_0)+2\,f_{n-1}(2n-1-k_0,2n-1-m_0)
=0.\quad\cr
}
$$
But by $(R\,3)$ written at
$(m,k)=(2n-1-k_0,2n+1-m_0)$ we have:
$$\displaylines{
f_{n}(2n-1-k_0,2n+1-m_0)-
2\,f_{n}(2n-k_0,2n+1-m_0)\hfill\cr
\hfill{}-f_{n}(2n-k_0+1,2n+1-m_0)
+2\,f_{n-1}(2n-1-k_0,2n-1-m_0)=0.\cr}
$$
By comparing the last two equations we
conclude that
$$
f_{n}(m_0,k_0)=f_{n}(2n-k_{0}+1,2n+1-m_0),
$$
which means that (1.12) now holds for $(m,k)=(m_0,k_0)$.

\medskip
For the entries of~$M_{n}$ lying {\it below the diagonal} we
proceed in the same manner and adopt the notation:
$$\leqalignno{{\rm Row}_{2n+1-i}
&=\{(2n+1-i,i),\ldots,(2n+1-i,2n-i)\};\cr
{\rm Col}_{i}
&=\{(i+1,i),(i+2,i),\ldots,(2n+1-i,i)\};\cr
}
$$
for $i=1,2,\ldots,n-1$.
Again, ${\rm Row}_{2n+1-i}$ and ${\rm Col}_{i}$ have the cell
$(2n+1-i,i)$ in common.

Let $(m_0,k_0)$ belong to ${\rm Col}_j$ for some~$j$ such
that $1\le j\le n-1$. Further, assume that (1.12) holds for all 
$(m,k)\in{\rm Col}_{1}\cup\cdots\cup{\rm Col}_{j-1}$
and all $(m,k)\in {\rm Col}_j$ lying below
$(m_0,k_0)$, not including $(m_0,k_0)$, that is, $m> m_0$ and 
$k=k_0$. By symmetry, (1.12) also holds for all $(m,k)\in {\rm
Row}_{2n}\cup\cdots\cup{\rm Row}_{2n-j+2}$ and all $(m,k)\in
{\rm Row}_{2n+1-j}$ lying to the left of
$(2n+1-k_0,2n+1-m_0)$ not including the latter pair.

Now, the following relations hold:
$$
\displaylines{
f_{n}(m_0,k_0)\!=\!2f_{n}(m_0+1,k_0)
-f_{n}(m_0+2,k_0)
-2\,f_{n-1}(m_0,k_0),\cr
\noalign{\rightline{[by $(R\,1)$],}} 
f_{n-1}(m_0,k_0)=f_{n-1}(2n-1-k_0,2n-1-m_0),\cr
\noalign{\rightline{[by induction on~$n$],}}
f_{n}(m_0+1,k_0)=f_{n}(2n+1-k_0,2n-m_0),\cr
f_{n}(m_0+2,k_0)=f_{n}(2n+1-k_0,2n-1-m_0),\cr
\noalign{\rightline{[by the local induction].}}
\noalign{\hbox{Therefore,}}
\quad
f_{n}(m_0,k_0)=2\,f_{n}(2n+1-k_0,2n-m_0)
\hfill\cr
\hfill{}-f_{n}(2n+1-k_0,2n-1-m_0)
-2\,f_{n-1}(2n-1-k_0,2n-1-m_0).\quad\cr
}
$$
But by $(R\,4)$ written at
$(m,k)=(2n+1-k_0,2n-1-m_0)$ we have:
$$\displaylines{\noalign{\vskip-7pt}
f_{n}(2n+1-k_0,2n-1-m_0)=
2\,f_{n}(2n+1-k_0,2n-m_0)\hfill\cr
\hfill{}-f_{n}(2n+1-k_0,2n+1-m_0)
-2\,f_{n-1}(2n-1-k_0,2n-1-m_0).\cr
\noalign{\hbox{By comparing the last two equations we
conclude that}}
f_{2n+1}(m_0,k_0)=f_{2n+1}(2n+1-k_0,2n+1-m_0)),\cr
\noalign{\vskip-7pt}}
$$
which means that (1.12) now holds for $(m,k)=(m_0,k_0)$.\cqfd

\medskip
Define
$$
\Eoc(t):=\eoc(t),\quad{\rm but}\quad\Pom(t):=2n+1-\pom(t).
\leqno(6.2)
$$ 

\proclaim Theorem 6.2. Let ${\goth T}_{2n+1}(x,y):=\kern-8pt 
\sum\limits_{t\in {\goth T}_{2n+1}}
x^{\Eoc(t)}y^{\Pom(t)}$ be the generating polynomial for
the set ${\goth T}_{2n+1}$ by the pair of statistics
$(\Eoc,\Pom)$. Then,
$$
{\goth T}_{2n+1}(x,y)={\goth T}_{2n+1}(y,x).
\leqno(6.3)
$$

\proof This is a simple consequence of Theorem 1.2:
let $g_{n}(m,k):=\#\{\Eoc=m,\Pom=k\}$. Then,
$$
\eqalignno{
g_{n}(m,k)&=\#\{\eoc=m,\pom=2n+1-k\}=f_{n}(m,2n+1-k)\cr
&=f_{n}(k,2n+1-m)&\hbox{[by (1.12)]}\cr
&=g_{n}(k,m).\cr
\noalign{\hbox{Thus,}}
{\goth T}_{2n+1}(x,y)&=\sum_{t\in {\goth T}_{2n+1}}
x^{\Eoc(t)}y^{\Pom(t)}=\sum_{m,k}g_{n}(m,k)x^my^k\cr
&=\sum_{m,k}g_{n}(k,m)x^my^k={\goth T}_{2n+1}(y,x).\qed\cr}
$$

\medskip
\centerline{\bf 7. Further properties}

\medskip
Other properties of the Delta Sequence can be obtained by having a further look at the geometry of the strictly ordered binary trees. The sub- and superdiagonals of the matrices~$M_{n}$ for $n=2,3,4,5$ are equal, as can be seen in Fig.~1.2. For an arbitrary $n\ge 2$ we have the following.

\proclaim Property 7.1. Sub- and super diagonals are
equal:
$$
f_{n}(k+1,k)=f_{n}(k,k+1)\qquad (n\ge 2;\;1\le k\le 2n-1).
\leqno(7.1)
$$

\proof
First, note that $k$ and $(k+1)$ can be siblings in a tree from
${\goth T}_{2n+1,k,k+1}$, but never in a tree from
${\goth T}_{2n+1,k+1,k}$. Second, $k$ can be parent of $(k+1)$ in a tree from the latter set, but never in a tree from the former one.
Also, $f_{n}(2,1)=f_{n}(1,2)=0$ for $n\ge 2$ and for $k\ge 2$ we have the decompositions:
$$\leqalignno{
f_{n}(k,k+1)&=\quad
\arbrebb{}{}{\raise4pt\hbox{$\scriptstyle
2n+1$}}{k+1}{\raise5pt\hbox{$\rond$}}{\lower2pt\hbox{$\scriptstyle
$}}{}{\raise5pt\hbox{$\scriptstyle\bullet\ k$}}
+\qquad[\qquad
\arbrebb{\lower0pt\hbox{$\rond$}}{}{\raise4pt\hbox{$\scriptstyle
2n+1$}}{k+1}{}{\lower2pt\hbox{$\scriptstyle
$}}{}{\raise5pt\hbox{$\carre$}}\quad,\ \bullet\ k\ ];
\cr
\noalign{\medskip}
f_{n}(k+1,k)&=\quad
\arbrebb{\lower2pt\hbox{$\scriptstyle\; \bullet$}}{}{\raise4pt\hbox{$\scriptstyle
2n+1$}}{k}{\lower2pt\hbox{$\scriptstyle
\;k+1$}}{\lower2pt\hbox{$\scriptstyle
$}}{}{\raise5pt\hbox{$\rond$}}
+\qquad[\qquad
\arbrebb{\lower0pt\hbox{$\rond$}}{}{\raise4pt\hbox{$\scriptstyle
2n+1$}}{k}{}{\lower2pt\hbox{$\scriptstyle
$}}{}{\raise5pt\hbox{$\carre$}}\quad,\ \bullet\ (k+1)\ ].
\cr
}
$$
The first terms in the previous two equations are in bijection,
as well as the second ones, the notation ``$\bullet \ k$'' meaning that $k$ is the end of the minimal chain, following our convention on Tree Calculus (\cf. Section 3).\qed

\proclaim Property 7.2. We have the crossing equalities:
$$\displaylines{(7.2)\quad
f_{n}(k+1,k-1)+
f_{n}(k-1,k+1)\hfill\cr
\kern3cm{}=f_{n}(k+1,k)+
f_{n}(k-1,k)\hfill\cr
\kern3cm{}=f_{n}(k,k+1)+f_{n}(k,k-1)
\quad (2\le k\le 2n-1).\hfill\cr}
$$

\goodbreak
Reporting to Fig. 1.2 the involved entries in the first identities are located on the
four bullets drawn in the following diagram.
$$
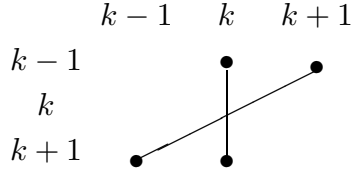
\vbox{\vskip1.3cm\offinterlineskip
\centerput(0,10){$k-1$}
\centerput(10,10){$k$}
\centerput(20,10){$k+1$}
\centerput(10,5){$\bullet$}
\centerput(20,4.5){$\bullet$}
\segment(20,5)\dir(-2,-1)\long{20}
\centerput(0,-6){$\bullet$}
\centerput(10,-6){$\bullet$}
\segment(10,5)\dir(0,-1)\long{10}
\centerput(-10,5){$k-1$}
\centerput(-10,0){$k$}
\centerput(-10,-5){$k+1$}
\vskip.7cm
}\hskip1cm
$$

\proof
Let $i$, $j$ be two different integers from the set
$\{(k-1), k,(k+1)\}$. Say that $i$ and~$j$ are {\it connected} in a tree~$t$, if the tree contains the edge $i$---$j$, or 
if $i$ and $j$ are brothers and
one of them is the end of the minimal chain of~$t$. 
Each of the four ingredients of the previous identity 
is now decomposed into five terms, depending on whether the nodes $(k-1)$, $k$, $(k+1)$ are connected or not, namely: no connectedness; only $k, (k+1)$ connected; $(k-1), k$ connected; $(k-1),(k+1)$ connected; all connected. Thus,
$$
\leqalignno{
f_{n}(k-1,k)&=
[\quad\arbrea{k-1}{\lower0pt\hbox{$\rond$}}{}{\lower2pt
\hbox{$\bullet$}}\ ,\quad
\arbrea{}{\lower0pt\hbox{$\carre$}}{k}{\raise2pt
\hbox{$\scriptstyle 2n+1$}}
,k+1]
+[\quad\arbrea{k-1}{\lower0pt\hbox{$\rond$}}{}{\lower2pt
\hbox{$\bullet$}}\ ,\qquad
\arbrea{2n+1}{\raise1pt\hbox{$\carre$}}{k}{\lower2pt
\hbox{}}\kern-5pt\raise12pt\hbox{$\scriptstyle
k+1$}\quad]\cr 
\noalign{\medskip}
&\kern-10pt{}+
[\quad
\arbrebb{\lower0pt\hbox{$\rond$}}{}{\raise4pt\hbox{$\scriptstyle
2n+1$}}{k}{}{}{}{\raise5pt\hbox{$\scriptstyle\bullet\  k-1$}}\ ,\ k+1\ ]
+[\quad\arbrea{k-1}{\lower4pt\hbox{$\rond\ \atop k+1$}}{}{\lower2pt
\hbox{$\bullet$}}\ ,\quad
\arbrea{}{\lower0pt\hbox{$\carre$}}{k}{\raise2pt
\hbox{$\scriptstyle 2n+1$}}]
+\qquad
\arbrebb{}{}{\raise4pt\hbox{$\scriptstyle
2n+1$}}{k}{\hbox{$\rond\atop k+1$}\hskip-5pt}{}{}{\raise5pt\hbox{$\scriptstyle\bullet\ k-1$}}
\cr
&:=A_1+A_2+A_3+A_4+A_5;\cr}
$$

\goodbreak
$$
\leqalignno{
f_{n}(k+1,k)&=
[\ \arbrea{k+1}{\lower0pt\hbox{$\rond$}}{}{\lower2pt
\hbox{$\bullet$}}\ ,k-1,\arbrea{2n+1}{\lower0pt\hbox{$\rond$}}{k}{\lower2pt
\hbox{$\bullet$}}]\quad
+\quad[\
\arbrea{k+1}{\lower-5pt\hbox{$\kern-15pt\textstyle\  \atop
2n+1$}}{k}{\lower2pt
\hbox{$\bullet$}}\
,k-1]\cr
\noalign{\medskip}
&\kern-1.5cm{}
+[\ \arbrea{k+1}{\lower0pt\hbox{$\rond$}}{}{\lower2pt
\hbox{$\bullet$}}\
,\qquad\quad \arbrebb{\lower2pt\hbox{}}{\rond}%
{}{k}{2n+1}{\kern-3pt \carre}{k-1}{}\quad]
+[\
\arbrea{k+1}{\lower0pt\hbox{$\carre$}}{k-1}{\lower2pt
\hbox{$\bullet$}}\
,\qquad
\arbrea{2n+1}{\lower0pt\hbox{$\rond$}}{k}{\lower2pt
\hbox{}}]\quad
+\quad\arbrebb{\lower2pt\hbox{$\;\bullet$}}{\rond}%
{}{k}{k+1}{\kern-3pt{\textstyle} 2n+1}{k-1}{}
\cr
\noalign{\smallskip}
&\kern-1cm{}:=B_1+B_2+B_3+B_4+B_5;\cr
\noalign{\smallskip}
f_{n}(k+1,k-1)&=
[\ \arbrea{k+1}{\lower0pt\hbox{$\rond$}}{}{\lower2pt
\hbox{$\bullet$}},\quad
\arbrea{}{\lower0pt\hbox{$\carre$}}{k-1}{\raise2pt
\hbox{$\scriptstyle 2n+1$}},\ k]
+\qquad[\ \arbrea{k+1}{\carre}{k}{\lower2pt
\hbox{$\bullet$}}\quad
,\ \arbrea{}{\hskip4pt
\hbox{$\scriptstyle
2n+1$}}{k-1}{\rond}\quad]\cr
\noalign{\medskip}
&\kern-1cm{}+
[\ \arbrea{k+1}{\lower0pt\hbox{$\rond$}}{}{\lower2pt
\hbox{$\bullet$}}\ ,\qquad
\arbrea{2n+1}{\raise 1.5pt\hbox{$\carre$}}{}{\lower2pt
\hbox{}}\kern-13pt \raise4pt\hbox{$\scriptstyle k-1$}
\kern-7pt \raise 11pt\hbox{$\scriptstyle k$}
\quad ]+
[\
\arbrea{k+1}{\lower-5pt\hbox{$\textstyle\  \atop
2n+1$}}{k-1}{\lower2pt
\hbox{$\bullet$}}\
\quad ,k]
+\quad\arbrebb{\lower2pt\hbox{$\;\bullet$}}{\lower
2.5pt\hbox{$\ $} }{\rond}{k}{k+1}{}{k-1}{2n+1}\cr
\noalign{\smallskip}
&\kern-1cm{}:=C_1+C_2+C_3+C_4+C_5;\cr
\noalign{\medskip}
f_{n}(k-1,k+1)&= [\
\arbrea{k-1}{\lower0pt\hbox{$\rond$}}{}{\lower2pt
\hbox{$\bullet$}}\ ,\ k\ ,\arbrea{}{\lower0pt\hbox{$\carre$}}{k+1}{\raise2pt
\hbox{$\scriptstyle 2n+1$}}\quad ]
+[\
\arbrea{k-1}{\lower0pt\hbox{$\carre$}}{}{\lower2pt
\hbox{$\bullet$}},\qquad
\arbrebb{\triang}{\rond }{}{k+1}{}{2n+1}{k}{}\ ]
\cr
\noalign{\medskip}
&\kern-1cm{}+\ 
[\ \arbrea{k-1}{\lower4pt\hbox{$\carre\atop k$}}{}{\lower2pt
\hbox{$\bullet$}}\quad
,\ \arbrea{}{\raise2pt
\hbox{$\scriptstyle
2n+1$}}{k+1}{\rond}\quad]
+[\ \quad\arbrebb{\rond}{\lower
2.5pt\hbox{$\ $} }{\raise3pt\hbox{$\scriptstyle 2n+1$}}{k+1}{}{}{}{\hskip-10pt\hbox{$\scriptstyle\bullet\atop \hskip5pt  k-1$}}\quad
,k\ ]
+\quad\arbrebbc{2n+1}{\carre}{}{k+1}{}{\hskip-3pt\hbox{$\rond$} }{k}{}
{\hskip14pt\lower4pt\hbox{$\scriptstyle \bullet\ k-1$}}{}\quad\cr
\noalign{\smallskip}
&\kern-1cm:=D_1+D_2+ D_3+ D_4+ D_5.\cr
}
$$

Now, the following identities hold: $A_1= C_1$, 
$A_2= C_3$, $A_3= D_{4}$,
$A_4= D_{3}$, $2\,A_5= D_{5}$; $B_1= D_1$, $B_3= D_2$,
$B_5= C_5$. Moreover,
$$
\leqalignno{B_{4}-C_{2}&=
\quad\lower9pt\hbox{$\arbrebbz{\hskip-15pt 
\lower 2pt\hbox{$\scriptstyle k+1\ \bullet$}}{\quad k}{}{k-1}{\lower2pt\hbox{$\carre$}}{ \; 2n+1}{}
{}{\ \rond}$}\quad=D_{5}=2\,A_{5;}
\cr
C_{4}-B_{2}&=\qquad\arbrebb{\lower2pt\hbox{$\;\bullet$}}{\lower
2.5pt\hbox{$\ $} }{\raise2pt\hbox{$\scriptstyle 2n+1$}}{k-1}{k+1}{}{}{%
\kern -3pt\hbox{\raise0pt\hbox{$\rond\atop k$}}
}
\qquad=A_{5}.\cr
}
$$
Altogether, $\sum_{i}(A_{i}+B_{i})-\sum_{i}(C_{i}+D_{i})
=(B_{4}-C_{2})-(C_{4}-B_{2})+(A_{5}-D_{5})
=2\,A_{5}-A_{5}+(A_{5}-2\,A_{5})=0$.\qed

\goodbreak

\proclaim Property 7.3. The row sums $f_{n}(m,\brullet)$ form a Poupard Triangle, the initial conditions being:
$f_{0}(1,\brullet)=1$, $f_{n}(1,\brullet)=0$ and 
$f_{n}(2,\brullet)=2\,\sum_{m}f_{n-1}(m,\brullet)$ $(n\ge 1)$;
and the finite difference system:
$$
\Deltaa_m{\kern-3pt}^2 f_{n}(m,\brullet)+2\,f_{n-1}(m,\brullet)=0
\quad 
(1\le m\le 2n-1).\leqno(7.3)
$$
The column sums $f_{n}(\brullet,k)$ form a Poupard Triangle, the initial conditions being:
$f_{0}(\brullet,0)=1$, $f_{n}(\brullet,0)=0$ and 
$f_{n}(\brullet,1)=2\,\sum_{k}f_{n-1}(\brullet,k)$ $(n\ge 1)$;
and the finite difference system:
$$
\Deltaa_k{\kern-3pt}^2 f_{n}(\brullet,k)+2\,f_{n-1}(\brullet,k)=0
\quad 
(0\le k\le 2n-2).\leqno(7.4)
$$

There are several proofs of this Property. First, the methods developed in Section~3 can be readapted by disregarding the conditions involving the $\pom$-statistic. Here, we simply
work out a specialization of the recurrence relations $(R\,1)$---$(R\,4)$, that makes use of the previous two properties. Besides, we only prove the first part of the property that deals with the row sums.

\medskip
\proof
For $1\le m\le 2n-2$ we have:
$$\eqalignno{
\Deltaa_m{\kern-3pt}^2 f_{n}(m,\brullet)
&=\sum_{1\le k\le 2n}\Deltaa_m{\kern-3pt}^2f_{n}(m,k)\cr
&=\sum_{1\le k\le m-1}\Deltaa_m{\kern-3pt}^2 f_{n}(m,k)\cr
&\kern.8cm{}+\Deltaa_m{\kern-3pt}^2 (f_{n}(m,m)+f_{n}(m,m+1)
+f_{n}(m,m+2))\cr
&\kern.8cm{}+\sum_{m+3\le k\le 2n}
\Deltaa_m{\kern-3pt}^2f_{n}(m,k) \cr
&=\sum_{1\le k\le m-1}-2f_{n-1}(m,k)\cr
&\kern.8cm{}+f_{n}(m,m)-2f_{n}(m+1,m)+f_{n}(m+2,m)\cr
&\kern.8cm{}+f_{n}(m,m+1)-2f_{n}(m+1,m+1)+f_{n}(m+2,m+1)\cr
&\kern.8cm{}+f_{n}(m,m+2)-2f_{n}(m+1,m+2)+f_{n}(m+2,m+2)\cr
&\kern.8cm{}+\sum_{m+3\le k\le 2n}-2f_{n-1}(m,k-2).\cr
}$$
In the previous sum the diagonal terms vanish. Also,
$f_{n}(m+1,m)=f_{n}(m,m+1)$, $f_{n}(m+2,m+1)=f_{n}(m+1,m+2)$ by
(7.1). The sum of the nine intermediate terms becomes:
$-f_{n}(m+1,m)+f_{n}(m+2,m)-f_{n}(m+1,m+2)+f_{n}(m,m+2)$, which is~0 by (7.2). Hence,
$$\eqalignno{
\Deltaa_m{\kern-3pt}^2 f_{n}(m,\brullet)
&=-2\sum_{1\le k\le m-1}f_{n-1}(m,k)
-2\sum_{m+1\le k\le 2n-2}f_{n-1}(m,k)\cr
&=-2f_{n-1}(m,\brullet).\cr}
$$

\goodbreak
On the other hand, by $(I\,2)$ and (1.11),
$$\eqalignno{
\quad f_{n}(2n-1,\brullet)
&=\sum\limits_{1\le k\le 2n}f_{n}(2n-1,k)=\sum\limits_{1\le k\le 2n-2}
(f_{n-1}(k,\brullet)+f_{n-1}(\brullet,k))\cr
&=2\,f_{n-1}(\brullet,\brullet):=2\,\#{\goth T}_{2n-1}
=2\,T_{2n-1}/2^{n-1};\cr
f_{n}(2n,\brullet)
&=\sum\limits_{1\le k\le 2n}f_{n}(2n,k)=\sum\limits_{1\le k\le 2n-2}
f_{n-1}(k,\brullet)\cr
&=f_{n-1}(\brullet,\brullet)=\#{\goth T}_{2n-1}
=T_{2n-1}/2^{n-1};\cr
f_{n}(2n+1,\brullet)&=0.\cr
}
$$
The last three evaluations imply that $\Deltaa_m{\kern-3pt}^2 f_{n}(2n-1,\brullet)=-2\,f_{n-1}(2n-1,\brullet)=0$.\cqfd

\bigskip
\centerline{\bf 8. Other equivalent definitions for the Delta sequence}

\medskip
Definitions 1.1 and 1.2 have been shown to be equivalent to characterize a Delta Sequence. Other combinations of the recurrence relations $(R\,1)$--$(R\,4)$, together with the initial conditions $(I\,1)$--$(I\,4)$, can be used. We describe them by means of squares and arrows, as was done in Fig.~1.3:
$$
\vbox{\vskip2.2cm\offinterlineskip
\segment(0,0)\dir(1,0)\long{16}
\segment(0,2)\dir(1,0)\long{16}
\segment(0,0)\dir(0,1)\long{16}
\segment(16,0)\dir(0,1)\long{16}
\segment(0,14)\dir(1,0)\long{16}
\segment(0,16)\dir(1,0)\long{16}
\segment(0,16)\dir(1,-1)\long{16}
\centerput(8,13){$\downarrow$}
\centerput(8,15.5){$\uparrow$}
\centerput(8,-1){$\downarrow$}
\centerput(8,1.5){$\uparrow$}
\centerput(-6,3){$\scriptstyle {\rm Row}_{2n-1}$}
\centerput(-4.5,0){$\scriptstyle {\rm Row}_{2n}$}
\centerput(-4,15){$\scriptstyle {\rm Row}_{1}$}
\centerput(-4,13){$\scriptstyle {\rm Row}_{2}$}
\centerput(8,-4){$\scriptstyle (R\,1)$}
\centerput(8,20){$\scriptstyle (R\,3)$}
\centerput(8,-8){\eightrm Definition 1.3}
\vskip .7cm
}
\hskip3.5cm
\vbox{\vskip2cm\offinterlineskip
\segment(0,0)\dir(1,0)\long{16}
\segment(0,0)\dir(0,1)\long{16}
\segment(16,0)\dir(0,1)\long{16}
\segment(14,0)\dir(0,1)\long{16}
\segment(0,16)\dir(1,0)\long{16}
\segment(2,0)\dir(0,1)\long{16}
\segment(0,16)\dir(1,-1)\long{16}
\centerput(1,8){$\longleftrightarrow$}
\centerput(22,8){$\scriptstyle (R\,2)$}
\centerput(15,8){$\longleftrightarrow$}
\centerput(-7,8){$\scriptstyle (R\,4)$}
\centerput(14,17){$\scriptstyle {\rm Col}_{2n-1}$}
\centerput(20,14.5){$\scriptstyle {\rm Col}_{2n}$}
\centerput(-2,17){$\scriptstyle {\rm Col}_{1}$}
\centerput(4,17){$\scriptstyle {\rm Col}_{2}$}
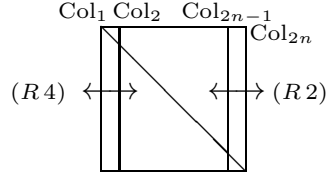
\centerput(8,-8){\eightrm Definition 1.4}
\vskip .7cm
}\hskip1cm$$

Other initial conditions than $(I\,1)$---$(I\,4)$ can be introduced. They will be denoted by 
$(SW\,1,4)$, $(N\!E\,2,3)$, as they refer only to the South-West and North-East corners of the matrices: 

\vskip-8pt
{\eightpoint
$$\leqalignno{\qquad
\petitematrice{f_{n}(2n-1,1)\ &\ f_{n}(2n-1,2)\cr
f_{n}(2n,1)&f_{n}(2n,2)\cr}
&\!=\!\petitematrice{f_{n-1}(\brullet,1)\ &\ f_{n-1}(2,\brullet)
+f_{n-1}(\brullet,2)\cr
0&f_{n-1}(\brullet,1)\cr};&(SW)\cr
\noalign{\medskip}
\petitematrice{f_{n}(1,2n-1)\ &\ f_{n}(1,2n)\cr
f_{n}(2,2n-1)&f_{n}(2,n)\cr}&=
\petitematrice{0&\ 0\cr
f_{n-1}(2,\brullet)&\ 0\cr}
=\petitematrice{0&\ 0\cr
f_{n-1}(2n-2,\brullet)&\ 0\cr}.&(NE)\cr}
$$

}
When one of those two conditions $(SW)$, $(N\!E)$ is involved, two recurrence relations among $(R\,1)$--$(R\,4)$ are needed to build up an equivalent definition. In Definition~1.5 for instance, 
$(R\,1)$ and $(R\,4)$ are to be associated with
$(SW)$. We then get five further equivalent definitions:

$$
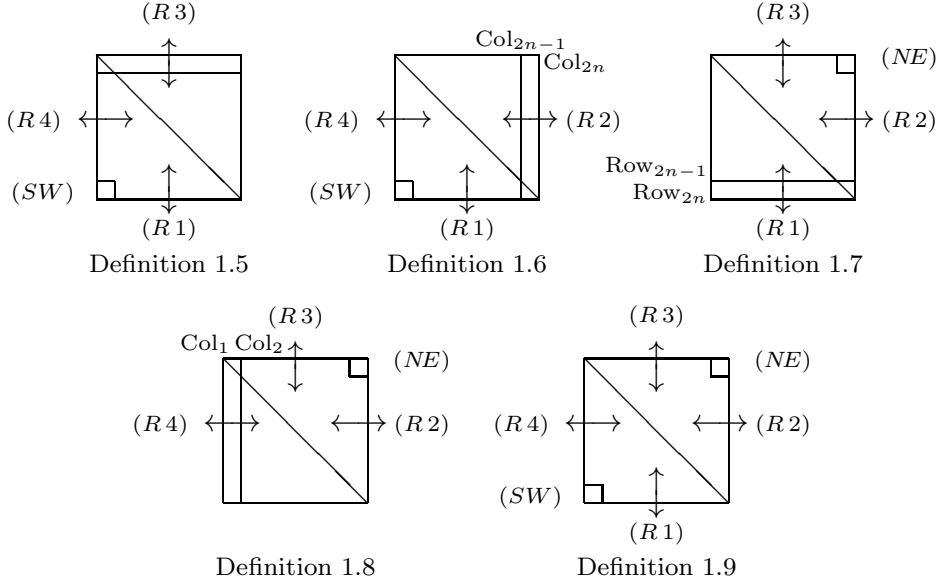
\vbox{\vskip2cm\offinterlineskip
\segment(0,0)\dir(1,0)\long{16}
\segment(0,2)\dir(1,0)\long{2}
\segment(0,0)\dir(0,1)\long{16}
\segment(16,0)\dir(0,1)\long{16}
\segment(0,14)\dir(1,0)\long{16}
\segment(0,16)\dir(1,0)\long{16}
\segment(2,0)\dir(0,1)\long{2}
\segment(0,16)\dir(1,-1)\long{16}
\centerput(1,8){$\longleftrightarrow$}
\centerput(8,13){$\downarrow$}
\centerput(8,15.5){$\uparrow$}
\centerput(8,-1){$\downarrow$}
\centerput(8,1.5){$\uparrow$}
\centerput(-6,0){$\scriptstyle (SW)$}
\centerput(-7,8){$\scriptstyle (R\,4)$}
\centerput(8,-4){$\scriptstyle (R\,1)$}
\centerput(8,20){$\scriptstyle (R\,3)$}
\centerput(8,-8){\eightrm Definition 1.5}
\vskip .7cm
}\hskip2cm
\hskip1.3cm
\vbox{\vskip2cm\offinterlineskip
\segment(0,0)\dir(1,0)\long{16}
\segment(0,2)\dir(1,0)\long{2}
\segment(0,0)\dir(0,1)\long{16}
\segment(16,0)\dir(0,1)\long{16}
\segment(14,0)\dir(0,1)\long{16}
\segment(0,16)\dir(1,0)\long{16}
\segment(2,0)\dir(0,1)\long{2}
\segment(0,16)\dir(1,-1)\long{16}
\centerput(1,8){$\longleftrightarrow$}
\centerput(8,-1){$\downarrow$}
\centerput(8,1.5){$\uparrow$}
\centerput(-6,0){$\scriptstyle (SW)$}
\centerput(-7,8){$\scriptstyle (R\,4)$}
\centerput(8,-4){$\scriptstyle (R\,1)$}
\centerput(22,8){$\scriptstyle (R\,2)$}
\centerput(15,8){$\longleftrightarrow$}
\centerput(14,17){$\scriptstyle {\rm Col}_{2n-1}$}
\centerput(20,14.5){$\scriptstyle {\rm Col}_{2n}$}

\centerput(8,-8){\eightrm Definition 1.6}
\vskip .7cm
}\hskip3.5cm
\vbox{\vskip2cm\offinterlineskip
\segment(0,0)\dir(1,0)\long{16}
\segment(0,2)\dir(1,0)\long{16}
\segment(0,0)\dir(0,1)\long{16}
\segment(16,0)\dir(0,1)\long{16}
\segment(14,14)\dir(0,1)\long{2}
\segment(14,14)\dir(1,0)\long{2}
\segment(0,16)\dir(1,0)\long{16}
\segment(0,16)\dir(1,-1)\long{16}
\centerput(8,13){$\downarrow$}
\centerput(8,15.5){$\uparrow$}
\centerput(8,-1){$\downarrow$}
\centerput(8,1.5){$\uparrow$}
\centerput(22,15){$\scriptstyle (N\!E)$}
\centerput(8,20){$\scriptstyle (R\,3)$}
\centerput(8,-4){$\scriptstyle (R\,1)$}
\centerput(22,8){$\scriptstyle (R\,2)$}
\centerput(15,8){$\longleftrightarrow$}
\centerput(-6,3){$\scriptstyle {\rm Row}_{2n-1}$}
\centerput(-4.5,0){$\scriptstyle {\rm Row}_{2n}$}
\centerput(8,-8){\eightrm Definition 1.7}
\vskip .7cm
}\hskip1.7cm
$$

$$
\vbox{\vskip2cm\offinterlineskip
\segment(0,0)\dir(1,0)\long{16}
\segment(2,0)\dir(0,1)\long{16}
\segment(0,0)\dir(0,1)\long{16}
\segment(16,0)\dir(0,1)\long{16}
\segment(14,14)\dir(0,1)\long{2}
\segment(14,14)\dir(1,0)\long{2}
\segment(0,16)\dir(1,0)\long{16}
\centerput(-7,8){$\scriptstyle (R\,4)$}
\segment(0,16)\dir(1,-1)\long{16}
\centerput(1,8){$\longleftrightarrow$}
\centerput(8,13){$\downarrow$}
\centerput(8,15.5){$\uparrow$}
\centerput(22,15){$\scriptstyle (N\!E)$}
\centerput(8,20){$\scriptstyle (R\,3)$}
\centerput(22,8){$\scriptstyle (R\,2)$}
\centerput(15,8){$\longleftrightarrow$}
\centerput(-2,17){$\scriptstyle {\rm Col}_{1}$}
\centerput(4,17){$\scriptstyle {\rm Col}_{2}$}
\centerput(8,-8){\eightrm Definition 1.8}
\vskip .7cm
}\hskip4cm
\vbox{\vskip2cm\offinterlineskip
\segment(0,0)\dir(1,0)\long{16}
\segment(0,2)\dir(1,0)\long{2}
\segment(0,0)\dir(0,1)\long{16}
\segment(16,0)\dir(0,1)\long{16}
\segment(0,16)\dir(1,0)\long{16}
\segment(2,0)\dir(0,1)\long{2}
\segment(0,16)\dir(1,-1)\long{16}
\centerput(1,8){$\longleftrightarrow$}
\centerput(8,-1){$\downarrow$}
\centerput(8,1.5){$\uparrow$}
\centerput(-6,0){$\scriptstyle (SW)$}
\centerput(-7,8){$\scriptstyle (R\,4)$}
\centerput(8,-4){$\scriptstyle (R\,1)$}
\centerput(22,8){$\scriptstyle (R\,2)$}
\centerput(15,8){$\longleftrightarrow$}
\centerput(8,-8){\eightrm Definition 1.9}
\segment(14,14)\dir(0,1)\long{2}
\segment(14,14)\dir(1,0)\long{2}
\centerput(8,13){$\downarrow$}
\centerput(8,15.5){$\uparrow$}
\centerput(22,15){$\scriptstyle (N\!E)$}
\centerput(8,20){$\scriptstyle (R\,3)$}
\vskip .7cm
}\hskip1.7cm
$$

We do not reproduce any proofs for those equivalences, but point out the fact that our Tree Calculus requires that each initial condition be combinatorially interpreted, as was done in Sections~2 and~5. 

\bigskip
\centerline{\bf 9. Generating functions for the Delta sequence}

\medskip
1. {\it Poupard matrices}.\quad
Let $G=(g_{i,j})$ $(i\ge 0,\,j\ge 0)$ be an infinite matrix with nonnegative integral entries. Say that~$G$ is a {\it Poupard matrix}, if for every $i\ge 0$, $j\ge 0$ the following identity holds:
$$
g_{i,j+2}-2\,g_{i+1,j+1}+g_{i+2,j}+2\,g_{i,j}=0.\leqno(9.1)
$$
Let $G(x,y):=\sum\limits_{i\ge 0,\,j\ge 0}g_{i,j}\,(x^i/i)\,(y^j/j!)$;\quad $R_{i}(y):=\sum\limits_{j\ge 0}g_{i,j}\,(y^j/j!)$ $(i\ge 0)$;\quad $C_{j}(x):=\sum\limits_{i\ge 0}g_{i,j}\,(x^i/i)$ $(j\ge 0)$ be the exponential generating functions for the matrix itself, its rows and columns, respectively.

\proclaim Proposition 9.1. The following four properties are equivalent.\hfil\break\indent
(i) $G=(g_{i,j})$ $(i\ge 0,\,j\ge 0)$ is a Poupard matrix;
\hfil\break\indent
(ii) $R''_i(y)-2\,R'_{i+1}(y)+R_{i+2}(y) +2\,R_{i}(y)=0$
for all $i\ge 0$;
\hfil\break\indent
(iii) $C''_{j}(x)-2\,C'_{j+1}(x)+C_{j+2}(x)+2\,C_{j}(x)$ for all $j\ge 0$;
\hfil\break\indent
(iv) $G(x,y)$ satisfies the partial differential equation:
$$
{\partial^2 G(x,y)\over \partial x^2}
-2\,{\partial^2 G(x,y)\over \partial x\,\partial y}
+{\partial^2 G(x,y)\over \partial y^2}
+2\,G(x,y)=0.\leqno(9.2)
$$

\proof
It suffices to write
$R''_i(y)-2\,R'_{i+1}(y)+R_{i+2}(y) +2\,R_{i}(y)
=\sum\limits_{j\ge 0}(g_{i,j+2}-2\,g_{i+1,j+1}
+2\,g_{i+2,j}-2\,g_{i,j})(y^j/j!)$, and, similarly,
$C''_{j}(x)-2\,C'_{j+1}(x)+C_{j+2}(x)+2\,C_{j}(x)
=\sum\limits_{i\ge 0}(g_{i+2,j}-2\,g_{i+1,j+1}
+2\,g_{i,j+2}-2\,g_{i,j})(x^i/i!)$ to obtain the equivalence between the first three properties. As for the last one, simply note that
$G(x,y)=\smash{\sum_{i\ge 0}R_{i}(y)\,x^i/i!}
=\sum\limits_{j\ge 0}C_{j}(x)\,y^j/j!$ and make the appropriate derivations.\hskip-6pt\qed

\proclaim Proposition 9.2. We have
$$
G(x,y)=A(x+y)\,\cos(\sqrt 2\,y)+B(x+y)\,
\sin(\sqrt 2\,y),\leqno(9.3)
$$
where $A(x)$ and $B(x)$ are two arbitrary series.

\proof
Let  $\xi:=x+y$, $\eta:=y$. Then,
$\partial G/\partial x=\partial G/\partial\xi$;\quad 
$\partial G/\partial y=\partial G/\partial\xi+\partial G/\partial \eta$;\qquad
$\partial^2 G/\partial x^2=\partial^2 G/\partial \xi^2$;\qquad
$\partial^2 G/\partial y^2=\partial^2 G/\partial \xi^2
+2\,\partial^2 G/\partial \xi\,\partial \eta+
\partial^2 G/\partial \eta^2$;\qquad
$\partial^2 G/\partial x\,\partial y=
\partial^2 G/\partial \xi^2
+\partial^2 G/\partial \xi\,\partial \eta$.

Thus, (9.3) can be rewritten as
$$
{\partial^2 G\over \partial \eta^2}+2\,G=0.\leqno(9.4)
$$
The ordinary differential equation $G''+2\,G=0$, whose characteristic polynomial is $r^2+2=0$, has a general solution of the form $A\,\cos(\sqrt 2\,\eta)+
B\,\sin(\sqrt 2\,\eta)$, so that the general solution of (9.2) is exactly given by (9.3).~\qed

\medskip
The exact expression for the generating function $G(x,y)$ can then be derived, 
if $A(x+y)$ and $B(x+y)$ can be obtained by an independent calculation,
as done in the sequel.

\bigskip
2. {\it A sequence of Poupard matrices for the lower triangles}.\quad
The entries $f_{n}(m,k)$ $(1\le k<m\le 2n)$ from the lower triangles in the matrices $M_{n}$ $(n\ge 1)$ are now recorded as entries $\lambda^{(p)}_{i,j}$ $(p\ge 0,\,i\ge 0,\,j\ge 0)$ of infinite matrices $\Lambda^{(p)}=(\lambda^{(p)}_{i,j})$ $(i\ge 0,\,j\ge 0)$  as follows. 

Define
$$
\lambda^{(p)}_{i,j}:=\cases{0,&if $i+j\equiv p$ mod 2;\cr
f_{n}(m,k),&if $i+j\equiv p+1$ mod 2;\cr}\leqno(9.5)
$$
with $k:=j+1$, $m:=i+j+2$, $2n:=p+(i+j)+1$.
The latter equation makes sense, as $i+j$ and~$p$ are of different parity. The mapping $(p,i,j)\mapsto (n,m,k)$ is one-to-one, the reverse mapping being: $p:=2n-m+1$, $i:=m-k-1$, $j:=k-1$. 
Thus, for $p\ge 0$, the matrix $\Lambda^{(2p+1)}$ reads:
\smallskip
{\eightpoint
\vskip-14pt
$$\bordermatrix{&0&1&2&3&4&5\cr
0&f_{p+1}(2,1)&0&f_{p+2}(4,3)&0&f_{p+3}(6,5)&0&\cdots\cr
1&0&f_{p+2}(4,2)&0&f_{p+3}(6,4)&0&\cdots\cr
2&f_{p+2}(4,1)&0&f_{p+3}(6,3)&0&\cdots\cr
3&0&f_{p+3}(6,2)&0&\cdots\cr
4&f_{p+3}(6,1)&0&\cdots\cr
5&0&\cdots\cr
6&\cdots\cr
};
$$
}
and for $p\ge 1$, 
the matrix $\Lambda^{(2p)}$ is equal to

\vskip-14pt
{\eightpoint
$$
\bordermatrix{&0&1&2&3&4&5\cr
0&0&f_{p+1}(3,2)&0&f_{p+2}(5,4)&0&f_{p+3}(7,6)\cr
1&f_{p+1}(3,1)&0&f_{p+2}(5,3)&0&f_{p+3}(7,5)&\cdots\cr
2&0&f_{p+2}(5,2)&0&f_{p+3}(7,4)&\cdots\cr
3&f_{p+2}(5,1)&0&f_{p+3}(7,3)&\cdots\cr
4&0&f_{p+3}(7,2)&\cdots\cr
5&f_{p+3}(7,1)&\cdots\cr
6&\cdots\cr
}.
$$

}
 
 \proclaim Proposition 9.3. Every matrix $\Lambda^{(p)}$
$(p\ge 0)$ is a Poupard matrix.

\proof
Using Definition (9.5) we have
$$\displaylines{\quad
\lambda^{(p)}_{i,j+2}-2\, \lambda^{(p)}_{i+1,j+1}
+\lambda^{(p)}_{i+2,j} +2\,\lambda^{(p)}_{i,j}\hfill\cr
\kern2cm{}
=f_{n+1}(m+2,k+2)-2\,f_{n+1}(m+2,k+1)\hfill\cr
\kern5cm{}+f_{n+1}(m+2,k)+2\,f_{n}(m,k)
\hfill\cr
\kern2cm{}
=\Deltaa_{k}f_{n+1}(m+2,k)+2\,f_{n}(m,k)=0,\hfill\cr
}$$
by rule $(R\,4)$.\qed

\bigskip
3. {\it The first matrix $\Lambda^{(1)}$}.\quad
The counter-diagonal of $\Lambda^{(1)}$ at depth~$2i$ ($i\geq 0$) reads
$$\displaylines{
\lambda^{(1)}_{2i,0}=f_{i+1}(2i+2,1);\quad \lambda^{(1)}_{2i-1,1}=f_{i+1}(2i+2,2);\quad\cdots\hfill\cr
\hfill{}\cdots \ 
\lambda^{(1)}_{2i-j,j}=f_{i+1}(2i+2,j+1);\ \cdots\ 
\lambda^{(1)}_{0,2i}=f_{i+1}(2i+2,2i+1);\cr
}
$$
Those $(2i+1)$ terms are equal to the first $(2i+1)$ entries of the
$(2i+2)$-nd row of the matrix $M_{i+1}$, that is, 
to the single term $f_{1}(2,1)=1$ for $i=0$ and for $i\ge 1$ to
$$
f_{i}(1,\brullet)=0,\,f_{i}(2,\brullet),\,\ldots,\,
f_{i}(j+1,\brullet),\,\ldots,\,
f_{i}(2i,\brullet),\,0,
$$
by virtue of $(I\,2)$. Thus, $\Lambda^{(1)}$ is identical to the Poupard matrix $(f_{i}(j+1,\brullet))$ $(i,j\ge 0)$:

{\eightpoint
$$\displaylines{
\bordermatrix{&0&1&2&3&4&5&6\cr
0&1&0&0&0&0&0&0\cr
1&0&f_{1}(2,\brullet)&0&f_{2}(4,\brullet)&0&f_{3}(6,\brullet)\cr
2&0&0&f_{2}(3,\brullet)&0&f_{3}(5,\brullet)\cr
3&0&f_{2}(2,\brullet)&0&f_{3}(4,\brullet)\cr
4&0&0&f_{3}(3,\brullet)\cr
5&0&f_{3}(2,\brullet)\cr
6&0\cr
},\cr}
$$

\noindent
{\tenrm equal to}
$$
\displaylines{
\noalign{\vskip-10pt}
\bordermatrix{&0&1&2&3&4&5&6&\cdots\cr
0&1&0&0&0&0&0&0&\cdots\cr
1&0&1&0&1&0&4&\cdots\cr
2&0&0&2&0&8&\cdots\cr
3&0&1&0&10&\cdots\cr
4&0&0&8&\cdots\cr
5&0&4&\cdots\cr
6&0&\cdots\cr
\cdot&\cdots\cr
},\cr}
$$

}
\medskip 

\bigskip
4. {\it The first two columns of $\Lambda^{(p)}$}.\quad
In the sequel the column labeled~$j$ of $\Lambda^{(p)}$ will be denoted by $\Lambda^{(p)}_{\brullet,j}$ and the exponential generating function for that column by
$\Lambda^{(p)}_{\brullet,j}(x)=\sum_{i\geq 0} \lambda_{i,j}^{(p)} {x^i/ i!}$. 
Also,
$\Lambda^{(p)}(x,y):=\sum_{j\ge 0}\Lambda^{(p)}_{\brullet,j}(x)y^j/j!$ will be the double exponential generating function for the matrix $\Lambda^{(p)}$.

\proclaim Proposition 9.4. The first two columns of each matrix $\Lambda^{(p)}$ $(p\ge 2)$ are related to the columns of~$\Lambda^{(1)}$ by the identities:
$$
\Lambda^{(p)}_{\brullet,0}(x)
=\Lambda^{(1)}_{\brullet,p-1}(x),\qquad
\Lambda^{(p)}_{\brullet,1}(x)
={d\over dx}\Lambda^{(1)}_{\brullet,p-1}(x)
+\Lambda^{(1)}_{\brullet,p}(x).\leqno(9.6)
$$

\proof
For the first identity it suffices to prove 
$ \lambda_{i,0}^{(p)}=\lambda_{i, p-1}^{(1)}$, that is
$$
f_{n}(i+2, 1)=f_n(p+i+1,p)
$$
when $2n=p+i+1$. This is true by Theorem 1.2.
For the second identity it suffices to prove
$ \lambda_{i,1}^{(p)}=\lambda_{i+1, p-1}^{(1)}+
\lambda_{i, p}^{(1)}$, that is
$$
f_{n}(i+3, 2)=f_n(p+i+2,p)+f_n(p+i+2, p+1) \leqno (9.7)
$$
when $2n=p+i+2$. But by $(I\,4)$
$$
f_{n}(i+3, 2)=f_n(i+2,1)+f_n(i+3, 1),
$$
so that by Theorem 1.2, identity (9.7) holds.\qed
\goodbreak

5. {\it Solving  the partial differential equation}.\quad
By (9.6) 
$$\leqalignno{
\Lambda^{(p)}(x,y)\Bigm |_{\{y=0\}}
&=\Lambda^{(p)}_{\brullet,0}(x)=\Lambda^{(1)}_{\brullet,p-1}(x);\cr
{\partial \over \partial y} \Lambda^{(p)}(x,y)
\Bigm |_{\{y=0\}}
&=\Lambda^{(p)}_{\brullet,1}(x)
={d\over dx}\Lambda^{(1)}_{\brullet,p-1}(x)
+\Lambda^{(1)}_{\brullet,p}(x).\cr}
$$

As each $\Lambda^{(p)}$ is a Poupard matrix, we can use
identity (9.3), so that
$\Lambda^{(p)}(x,y)\Bigm |_{\{y=0\}}=
A(x)$ and
${\partial \Lambda^{(p)}(x,y)/\partial y}
=A'(x+y)\,\cos(\sqrt 2 y)-\sqrt 2\, A(x+y)\,\sin(\sqrt 2y)
+B'(x+y)\,\sin(\sqrt 2 y)+\sqrt 2\, B(x+y)\,\cos(\sqrt 2y)$. Hence,
$$\leqalignno{
	\Lambda^{(p)}(x,y)\Bigm |_{\{y=0\}}&=A(x)
=\Lambda^{(1)}_{\brullet,p-1}(x);\cr
{\partial \Lambda^{(p)}(x,y)\over\partial y}
\Bigm |_{\{y=0\}}&=A'(x)+\sqrt 2\,B(x)
={d\over dx}\Lambda^{(1)}_{\brullet,p-1}(x)
+\Lambda^{(1)}_{\brullet,p}(x).\cr}
$$
Consequently, $A(x)=\Lambda^{(1)}_{\brullet,p-1}(x)$ and
$B(x)=\Lambda^{(1)}_{\brullet,p}(x)/\sqrt 2$ and the general expression for $\Lambda^{p}(x,y)$ reads:
$$
\Lambda^{(p)}(x,y)
=\Lambda^{(1)}_{\brullet,p-1}(x+y)\,\cos(\sqrt 2\,y)+\Lambda^{(1)}_{\brullet,p}(x+y)\sin(\sqrt 2\,y)/\sqrt 2.\leqno(9.8)
$$

This expression still holds for $p=1$. We know that $\Lambda^{(1)}_{\brullet,0}(x)=1$. On the other hand, the coefficient of $x^{2k+1}/(2k+1)!$ $(k\ge 0)$ in $\Lambda^{(1)}_{\brullet,1}(x)$ is $f_{k+1}(2,\brullet)=T_{2k+1}/2^k$. Hence, $\Lambda^{(1)}_{\brullet,1}(x)=\sqrt 2\,\tan(x/\sqrt 2)$ and
$$
\leqalignno{
\Lambda^{(1)}(x,y)
&=\cos(\sqrt 2\,y)+\sqrt 2\,\tan\Bigl({x+y\over\sqrt 2}\Bigr){\sin(\sqrt 2\,y)\over \sqrt 2}\cr
&=\cos\Bigl({x-y\over \sqrt2}\Bigr)\Bigm/
\cos\Bigl({x+y\over \sqrt2}\Bigr),&(9.9)\cr
}
$$
a result already obtained by Poupard.
\goodbreak

{\it Remark}.\quad
For getting the solution for $\Lambda^{(1)}(x,y)$ we can also start with the general expression displayed in (9.3) and calculate $A$ and $B$ with the initial conditions
$\Lambda^{(1)}(0,y)=\Lambda^{(1)}(x,0)=1$. We find
$$\leqalignno{\Lambda^{(1)}(x,y)&=\cos(\sqrt 2 y)+{1-\cos(\sqrt 2(x+y))\over \sin(\sqrt 2(x+y))}\sin(\sqrt 2 y)\cr
&={\sin(\sqrt 2 x)+\sin(\sqrt 2 y)\over \sin(\sqrt 2(x+y))},&(9.10)\cr}
$$
an expression, which is naturally equal to the right-hand side of (9.9) (by a simple trigonometric calculation).

\medskip
We have not worked out other explicit formulas for $\Lambda^{(p)}(x,y)$ when $p\ge 3$, but only derived the exponential generating
function for those series, as explained in the next subsection.

\bigskip
5. {\it A generating function for the lower triangles}.\quad
By (9.9)
$$\leqalignno{
{\partial\over \partial y}
\Lambda^{(1)}(x,y)
&={\partial\over \partial y}
\cos\Bigl({x-y\over \sqrt2}\Bigr)\Bigm/
\cos\Bigl({x+y\over \sqrt2}\Bigr)\cr
&=
\sin(\sqrt 2x)\Bigm/ \sqrt 2\,\cos^2\Bigl(\displaystyle{x+y\over \sqrt2}\Bigr).&(9.11)\cr
\noalign{\vskip-5pt}}
$$

Let
$\displaystyle\Lambda(x,y,z):=\sum_{p\ge 1}
\Lambda^{(p)}(x,y){z^{p-1}\over (p-1)!}$.

\noindent
Then,
$$
\eqalignno{
\noalign{\vskip-5pt}
\Lambda(x,y,z)&
=\cos(\sqrt 2\,y)\sum_{p\ge 1}
{z^{p-1}\over (p-1)!}
\Lambda^{(1)}_{\brullet,p-1}(x+y)\cr
&{}\qquad+
{\sin(\sqrt 2\,y)\over\sqrt 2}
\sum_{p\ge 1}
{z^{p-1}\over (p-1)!}
\Lambda^{(1)}_{\brullet,p}(x+y)&\hbox{[by (9.8)]}\cr
&=\cos(\sqrt 2\,y)\,\Lambda^{(1)}(x+y,z)
+{\sin(\sqrt 2\,y)\over\sqrt 2}
{\partial\over \partial z}\Lambda^{(1)}(x+y,z),\cr
}
$$
since
$\Lambda^{(1)}(x,z)
\displaystyle=\sum_{p\ge 0}
{z^{p}\over p!}
\Lambda^{(1)}_{\brullet,p}(x)$. 

\goodbreak

By (9.11) we then get:
$$
\eqalignno{\Lambda(x,y,z)
&=\cos(\sqrt 2\,y)\,\cos\Bigl({x+y-z\over \sqrt2}\Bigr)\Bigm/
\cos\Bigl({x+y+z\over \sqrt2}\Bigr)\cr
&\kern3cm{}+\sin(\sqrt 2\,y) \sin(\sqrt2(x+y))
\Bigm/ 2\,\cos^2\Bigl(\displaystyle{x+y+z\over \sqrt2}\Bigr)\cr
&={\cos(\sqrt 2\, x)+\cos(\sqrt 2\, y)\,\cos(\sqrt 2\, z)\over
2\,\cos^2\Bigl(\displaystyle{x+y+z\over \sqrt2}\Bigr)}.\cr
}
$$
Now, express $\Lambda(x,y,z)$ as a series in the $f_{n}(m,k)$'s. By definition,
$$
\leqalignno{
\Lambda(x,y,z)&=\sum_{p,i,j}\lambda^{(p)}_{i,j}
{z^{p-1}\over (p-1)!}
{x^{i}\over i!}
{y^{j}\over j!}\quad(p\ge 1,\,i\ge 0,\,j\ge 0);\cr
\noalign{\hbox{so that by (9.5)}}
\Lambda(x,y,z)&=\sum_{k,m,n}
f_{n}(m,k){x^{m-k-1}\over (m-k-1)!}
{y^{k-1}\over (k-1)!}
{z^{2n-m}\over (2n-m)!},\cr
 }
$$
the latter sum over the set
$\{2\le k+1\le m\le 2n\}$. This achieves the proof of Theorem 1.4.

\medskip
Note that the right-hand side of (1.14) is symmetric in $y,z$, that is,
$\Lambda(x,y,z)=\Lambda(x,z,y)$. The change $y\leftrightarrow z$ in the left-hand side of (1.14) shows that
$$
f_{n}(2n+1-k,2n+1-m)=f_{n}(m,k),
$$
the symmetry proved for the entries 
$f_{n}(m,k)$ such that $m\ge k+1$.

\medskip
{\it Remark}.\quad
Let $z=0$ in (1.14). We get
$$\displaylines{
\sum_{2\le k+1\le 2n}
f_{n}(2n,k){x^{2n-k-1}\over (2n-k-1)!}
{y^{k-1}\over (k-1)!}
={\cos(\sqrt 2\, x)+\cos(\sqrt 2\, y)\over
2\,\cos^2\Bigl(\displaystyle{x+y\over \sqrt2}\Bigr)};\cr
\noalign{\hbox{or,}}
\noalign{\vskip-8pt}
\sum_{i\ge 0,\,j\ge 0}
\lambda^{(1)}_{{i,j}}{x^{i}\over i!}
{y^j\over j!}=\Lambda^{(1)}(x,y)
={\cos(\sqrt 2\, x)+\cos(\sqrt 2\, y)\over
2\,\cos^2\Bigl(\displaystyle{x+y\over \sqrt2}\Bigr)};\cr
}
$$
which is another expression for $\Lambda^{(1)}(x,y)$ than
(9.9) and (9.10).

\goodbreak
\bigskip
6. {\it A generating function for the upper triangles}.\quad
The entries $f_{n}(m,k)$ $(1\le m<k\le 2n)$ from the upper triangles in the matrices $M_{n}$ $(n\ge 1)$ are next recorded as entries $\omega^{(p)}_{i,j}$ $(p\ge 0,\,i\ge 0,\,j\ge 0)$ of infinite matrices $\Omega^{(p)}=(\omega^{(p)}_{i,j})$ $(i\ge 0, \,j\ge 0)$  as follows. 

Define
$$
\omega^{(p)}_{i,j}:=\cases{0,&if $i+j\not\equiv p$ mod 2;\cr
f_{n}(m,k),&if $i+j\equiv p$ mod 2;\cr}
$$
with $m:=p+1$, $k:=p+j+2$, $2n:=p+i+j+2$. Conversely, $i:=2n-k$, $j:=k-m-1$, $p:=m-1$. Thus, for $p\ge 0$

\smallskip
$\Omega^{(2p+1)}={}$
{\eightpoint
$$\leqalignno{\noalign{\vskip-10pt}
&\bordermatrix{&0&1&2&3\cr
0&0&\!f_{p+2}(2p+2,2p+4)&0&\cdots\cr
1&f_{p+2}(2p+2,2p+3)&0&f_{p+4}(2p+2,2p+5)&\cdots\cr
2&0&f_{p+3}(2p+2,2p+4)&\cdots\cr
3&f_{p+3}(2p+2,2p+3)&\cdots\cr
4&\cdots\cr
}\cr
}$$}

$\Omega^{(2p)}={}$
{\eightpoint
$$\leqalignno{\noalign{\vskip-10pt}
&\bordermatrix{&0&1&2&3\cr
0&f_{p+1}(2p+1,2p+2)&0&f_{p+2}(2p+1,2p+4)&0\cr
1&0&f_{p+2}(2p+1,2p+3)&0&\cdots\cr
2&f_{p+2}(2p+1,2p+2)&0&\cdots\cr
3&0&\cdots\cr
4&\cdots\cr
}\cr}
$$

}

{\it Remark}.\quad
The first rows of all matrices~$M_{n}$ are null, so that $\Omega^{(0)}$ is the infinite matrix with all entries equal to zero! 

\medskip
Also, write

$\Omega^{(1)}={}$

{\eightpoint
\vskip-10pt
$$\eqalignno{
&\bordermatrix{&0&1&2&3&4&5&6&7\cr
0&0&f_{2}(2,4)&0&f_{3}(2,6)&0&f_{4}(2,8)&0&f_{5}(2,10)\cr
1&f_{2}(2,3)&0&f_{3}(2,5)&0&f_{4}(2,7)&0&f_{5}(2,9)&\cdots\cr
2&0&f_{3}(2,4)&0&f_{4}(2,6)&0&f_{5}(2,8)&\cdots\cr
3&f_{3}(2,3)&0&f_{4}(2,5)&0&f_{5}(2,7)&\cdots\cr
4&0&f_{4}(2,4)&0&f_{5}(2,6)&\cdots\cr
5&f_{4}(2,3)&0&f_{5}(2,5)&\cdots\cr
6&0&f_{5}(2,4)&\cdots\cr
7&f_{5}(2,3)&\cdots\cr
}\cr
&\ =\bordermatrix{&0&1&2&3&4&5&6&7\cr
0&0&0&0&0&0&0&0&0\cr
1&1&0&1&0&4&0&34&\cdots\cr
2&0&2&0&8&0&68&\cdots\cr
3&1&0&10&0&94&\cdots\cr
4&0&8&0&104&\cdots\cr
5&4&0&94&\cdots\cr
6&0&68&\cdots\cr
7&34&\cdots\cr
}\cr
}
$$
}

\proclaim Proposition 9.6. Every matrix $\Omega^{(p)}$
$(p\ge 0)$ is a Poupard matrix.

Same proof as for Proposition 9.3.
\medskip

The row labeled~$i$ of $\Omega^{(p)}$ will be denoted by 
$\Omega^{(p)}_{i, \brullet}$ and the exponential generating function for that 
row by
$\Omega^{(p)}_{i, \brullet}(y)=\sum_{j\geq 0} \omega_{i,j}^{(p)} {y^j/ j!}$. 
Also,
$\Omega^{(p)}(x,y):=\sum_{i\ge 0}\Omega^{(p)}_{i,\brullet}(y)x^i/i!$ will be the double exponential generating function for the matrix $\Omega^{(p)}$.
As $x$ and $y$ play a symmetric role in (9.2), the solution in 
(9.2) may also be written
$$\eqalignno{
G(x,y)&=A(x+y)\,\cos(\sqrt 2\,x)+B(x+y)\,\sin(\sqrt 2\,x),\cr
\noalign{\hbox{
so that the generating function of each matrix 
$\Omega^{(p)}$ is of the form}}
\Omega^{(p)}(x,y)
&=A(x+y)\,\cos(\sqrt 2\,x)+B(x+y)\,\sin(\sqrt 2\,x).\cr}
$$
As the first row of each matrix $\Omega^{(p)}$ is the zero sequence, we have
$\Omega^{(p)}(0,y)=A(y)=0$.
Hence, $\Omega^{(p)}(x,y)=B(x+y)\,\sin(\sqrt 2\,x)$ and
$$
{\partial\over \partial x}\Omega^{(p)}(x,y)
=\Bigl({\partial\over \partial x}B(x+y)\Bigr)
\,\sin(\sqrt 2x)+\sqrt 2\ B(x+y)\,\cos(\sqrt 2x).$$
Therefore,
$$
\leqalignno{\noalign{\vskip-10pt}
\Omega^{(p)}_{1,\brullet}(y)&={\partial\over\partial x}
\Omega^{(p)}(x,y)|\Bigm|_{\{x=0\}}=\sqrt 2\,B(y)\cr
\noalign{\hbox{and then}}
\Omega^{(p)}(x,y)&={1\over \sqrt 2}\sin(\sqrt 2\,x)\,
\Omega^{(p)}_{1,\brullet}(x+y).&(9.12)\cr}
$$ 

The evaluation of $\Omega^{(1)}_{1,\brullet}(y)$ is easy, as the  row labeled~1 of $\Omega^{(1)}$ is $(1,0,1,0,4,0,34,0,\ldots\,)$, compared with $(0,1,0,1,0,4,0,34,0,\ldots\,)$, which is the sequence of the coefficients of the Taylor expansion of
$\sqrt 2\,\tan(y/\sqrt 2)$. 
In fact we have 
$\omega_{1,2j}^{(1)} = f_{j+2}(2, 2j+3) = f_{j+1}(2j+2, \brullet)
= T_{2j+1}/2^j.
$
Thus,
$$
\leqalignno{
\Omega^{(1)}_{1,\brullet}(y)&={d\over dy}
\,\sqrt 2\,\tan\Bigl({y\over \sqrt 2}\Bigr)
={1\over \cos^2(y/\sqrt 2)};\cr
\noalign{\hbox{so that}}
\Omega^{(1)}(x,y)
&={1\over \sqrt 2}\sin(\sqrt 2\,x)
{1\over \displaystyle\cos^2\Bigl({x+y\over\sqrt 2}\Bigr)}.&(9.13)\cr}
$$

\proclaim Proposition 9.7. For all $p\ge 1$ we have:
$$
\Omega^{(p)}_{1,\brullet}(y)
=\Omega^{(1)}_{p,\brullet}(y).\leqno(9.14)
$$

Same proof as for Proposition 9.4.
Now, define:
$$\Omega(x,y,z):=\displaystyle\sum_{p\ge 1}
\Omega^{(p)}(x,y){z^p\over p!}\leqno(9.15)
$$
and make use of (9.18)---(9.20):
$$
\leqalignno{
\Omega(x,y,z)
&={1\over \sqrt 2}\sin(\sqrt 2\,x)\sum_{p\ge 1}
\Omega^{(p)}_{1,\brullet}(x+y){z^p\over p!}\cr
&={1\over \sqrt 2}\sin(\sqrt 2\,x)\sum_{p\ge 1}
\Omega^{(1)}_{p,\brullet}(x+y){z^p\over p!}\cr
&={1\over \sqrt 2}\sin(\sqrt 2\,x)\,
\Omega^{(1)}(z,x+y)\cr
&=\sin(\sqrt 2\,x)\,
\sin(\sqrt 2\,z)\,
{1\over 2\,\displaystyle\cos^2\Bigl({x+y+z\over\sqrt 2}\Bigr)}.&\cr}
$$

As all the entries $\omega^{(p)}_{0,j}$ $(p\ge 1,\,j\ge 0)$ are null,
$$
\leqalignno{
\Omega(x,y,z)&=\sum_{p,i,j}\omega^{(p)}_{i,j}
{z^{p}\over p!}
{x^{i}\over i!}
{y^{j}\over j!}\quad(p\ge 1,\,i\ge 1,\,j\ge 0);\cr
\noalign{\hbox{so that by definition of the 
$\omega^{(p)}_{i,j}$'s}}
\Omega(x,y,z)&=\sum_{k,m,n}
f_{n}(m,k){x^{2n-k}\over (2n-k)!}
{y^{k-m-1}\over (k-m-1)!}
{z^{m-1}\over (m-1)!},&(9.16)\cr
 }
$$
the latter sum over the set
$\{2\le m+1\le k\le 2n-1\}$. This achieves the proof of Theorem 1.5.

\medskip

The right-hand side of (1.15) is symmetric in $x, z$, that is,
$\Omega(x,y,z)=\Omega(z,y,x)$. The change $x\leftrightarrow z$ in the left-hand side of (1.15) shows that
$$
f_{n}(2n+1-k,2n+1-m)=f_{n}(m,k),
$$
the symmetry proved for the entries 
$f_{n}(m,k)$ such that $m+1\le k$.

\vfill\eject
\vglue2cm

\centerline{\bf References}

{\eightpoint

\bigskip

\article  An1879|D\'esir\'e Andr\'e|D\'eveloppement de $\sec x$ et
$\tan x$|C. R. Math. Acad. Sci. Paris|88|1879|965--979|

\article An1881|D\'esir\'e Andr\'e|Sur les permutations
altern\'ees|J. Math. Pures et Appl.|7|1881|167--184|

\livre Co74|Comtet, Louis|\hskip-5pt Advanced 
Combinatorics|\hskip-5pt D.
Reidel/Dordrecht-Holland, Boston, {\oldstyle 1974}|

\divers FH12|Dominique Foata; Guo-Niu Han|Finite Difference Calculus
for Alternating Permutations, preprint, 15~p|

\divers Ha12|Guo-Niu Han|The Poupard Statistics on Tangent and Secant Trees, Strasbourg, preprint 12~p|

\livre Jo39|Charles Jordan|Calculus of Finite Differences|R\"ottig and
Romwalter,  Budapest, {\oldstyle 1939}|

\livre Ni23|Niels Nielsen|Trait\'e \'el\'ementaire des nombres
de Bernoulli|Paris, Gauthier-Villars, {\oldstyle 1923}|

\article Po89|Christiane Poupard|Deux propri\'et\'es des arbres
binaires ordonn\'es stricts|Europ. J. Combin.|10|1989|369--374|

\divers Sl07|N.J.A. Sloane|On-line Encyclopedia of Integer
Sequences,\hfil\break
{\tt
http://www.research.att.com/\char126 njass/sequences/}|

\divers Vi88|Xavier G. Viennot|S\'eries g\'en\'eratrices
\'enum\'eratives, chap.~3, Lecture Notes, 160~p., 1988, notes de
cours donn\'es
\`a l'\smash{\'E}cole Normale Sup\'erieure Ulm (Paris), UQAM (Montr\'eal,
Qu\'ebec) et Universit\'e de Wuhan (Chine)\hfil\break
{\tt
http://web.mac.com/xgviennot/Xavier\_Viennot/cours.html}|

\bigskip
\hbox{\vtop{\halign{#\hfil\cr
Dominique Foata \cr
Institut Lothaire\cr
1, rue Murner\cr
F-67000 Strasbourg, France\cr
\noalign{\smallskip}
{\tt foata@unistra.fr}\cr}}
\qquad
\vtop{\halign{#\hfil\cr
Guo-Niu Han\cr
I.R.M.A. UMR 7501\cr
Universit\'e de Strasbourg et CNRS\cr
7, rue Ren\'e-Descartes\cr
F-67084 Strasbourg, France\cr
\noalign{\smallskip}
{\tt guoniu.han@unistra.fr}\cr}}}

}

\bye